\documentclass[a4paper, final, 10pt]{article}
\usepackage{eurosym}
\usepackage{amsfonts,amsbsy,amssymb,amsthm}
\usepackage{amssymb}
\usepackage{amsthm}
\usepackage{graphicx}
\usepackage[format=plain, textfont=it]{caption}
\usepackage{subfig}
\usepackage{url}

\usepackage{pgfplotstable}
\usepackage{tikz}
\usetikzlibrary{patterns}
\usetikzlibrary{decorations.markings}

\tikzset{onearrow/.style={
        decoration={markings,
            mark= at position 1.0 with {\arrow{#1}} ,
        },
        postaction={decorate}
    }
}

\newcommand{\func}[1]{\mathop{\rm #1}}

\begin{document}

\title{ A new version of the
convexification method for a 1-D coefficient inverse problem with
experimental data\thanks{Supported by US Army Research Laboratory and US Army Research Office grant W911NF-15-1-0233 and by the Office of Naval Research grant N00014-15-1-2330.}}
\author{Michael V. Klibanov$^1$, Aleksandr E. Kolesov$^{1}$, Anders
Sullivan$^2$, Lam Nguyen $^2$}

\author{Michael V. Klibanov\thanks{The corresponding author} \thanks{Department of Mathematics \& Statistics, University of North Carolina at Charlotte, Charlotte, NC 28223, USA (mklibanv@uncc.edu, 
akolesov@uncc.edu)},  
Aleksandr E.\ Kolesov\footnotemark[2], Anders
Sullivan\thanks{U.S. Army Research Laboratory, Adelphi, MD
20783-1197 (lam.h.nguyen2.civ@mail.mil, anders.sullivan@us.army.mil)}, Lam Nguyen\footnotemark[3]}

\date{}

\maketitle

\begin{abstract}
A new version of the convexification method is developed analytically and
tested numerically for a 1-D coefficient inverse problem in the frequency
domain. Unlike the previous version, this one does not use the so-called
\textquotedblleft tail function", which is a complement of a certain
truncated integral with respect to the wave number. Globally strictly convex
cost functional is constructed with the Carleman Weight Function.
Global convergence of the gradient projection method to the correct solution
is proved. Numerical tests are conducted for both computationally simulated
and experimental data.
\end{abstract}



\section{Introduction}

\label{sec:1}

In this paper, we develop a new version of the so-called \emph{%
convexification} globally convergent numerical method for a 1-D
coefficient inverse problem.\ Next, we demonstrate its performance for both
computationally simulated and experimental data. The previous version of the
convexification was applied by this research group to the same experimental
data in \cite{KlibKol1}. Rather than collecting the data in a laboratory,
these data were collected in a realistic case of the cluttered environment
in the field. The data collection was performed by the Forward Looking Radar
of the US Army Research Laboratory \cite{Radar}. The goal of this radar is
to detect and identify flash explosive-like devices.

The radar community is currently relying only on the energy information of
radar images, see, e.g. \cite{Soumekh}. Unlike this, in the current paper,
so as in the previous one \cite{KlibKol1}, we compute estimates of
dielectric constants of targets. Our hope is that these estimates might help
in the future to develop new classification procedures, which would combine
the currently used energy information with the estimates of dielectric
constants. This combination, in turn might result in lower false alarm
rates. Our targets are three dimensional ones of course. On the other hand,
that radar can measure only one time dependent curve for each target. Thus,
what can be done at most by any data inversion technique is to estimate a
sort of an average of the dielectric constant for a target. We believe,
however, that even these estimates might be useful for the goal of
decreasing the false alarm rate. Thus, we model the wave propagation process
by the 1-D Helmholtz equation.

Convexification is the method, which constructs globally strictly convex
weighted Tikhonov-like functionals either for Coefficient Inverse Problems
(CIPs) or for ill-posed Cauchy problems for quasilinear PDEs. See our works \cite{KlibKol1, BK2, Klib95, Klib97, KT, KlibThanh, KlibKol2} for CIPs  and \cite{BakKlib, Klib15, KlibYag} for quasilinear PDEs. The key element of such a
functional is the presence of the Carleman Weight Function (CWF). This
is the weight function in the Carleman estimate for the principal part of
the corresponding differential operator.

Thus, convexification addresses the well known problem of multiple local
minima and ravines of conventional Tikhonov-like functionals for CIPs, see,
e.g. the paper \cite{Scales} for a numerical example of multiple local
minima. Convexification allows one to construct globally convergent
numerical methods for CIPs. We call a numerical method for a CIP \emph{%
globally }convergent if a theorem is proved, which claims that this method
delivers at least one point in a sufficiently small neighborhood of the
exact coefficient without any advanced knowledge of this neighborhood. Since
conventional least squares Tikhonov functionals are non convex, then they
usually have many local minima and ravines. This means that in order to
obtain the correct solution, using such a functional, one needs to start the
optimization process in a sufficiently small neighborhood of this solution,
i.e. one should work with a \emph{locally} convergent numerical method.
However, such a small neighborhood is rarely known in applications.

The convexification is a globally convergent numerical method, see Remark
4.1 in section 4. The idea of the convexification has roots in the method of
Carleman estimates for CIPs. This method was originated in the work of
Bukhgeim and Klibanov (1981) \cite{BukhKlib} as the tool of proofs of global
uniqueness and stability theorems for CIPs with the data resulting from a
single measurement event. The method of \cite{BukhKlib}  became quite
popular since then with many publications of many authors, see, e.g. a
survey in \cite{Ksurvey}, books \cite{BK1,BY,KT} and references cited
therein.

Another globally convergent numerical method for CIP with single measurement
data is the so-called \textquotedblleft tail functions" method, see, e.g. 
\cite{BK1,Kfreq}. This method is currently completely verified on
experimental data in a number of publications, see, e.g. \cite{Kfreq} and
reference cited therein.

The first publications on the convexification were \cite{Klib95,Klib97}.
Some theoretical gaps, which prevented one from the active numerical
studies, were fixed only recently in 
\cite{BakKlib}. As to the recent numerical results for the convexification
for CIPs, we refer to \cite{KlibKol1, KlibThanh} for the 1-D case and to 
\cite{KlibKol2} for the $3-$D case. Numerical results for the
convexification for ill-posed Cauchy problems for quasilinear PDEs were
published in \cite{BakKlib,KlibYag}. We also refer here to a recent
interesting version of the convexification for a CIP for the equation $%
u_{tt}=\Delta u+a\left( x\right) u$ in which one of initial conditions at $%
\left\{ t=0\right\} $ is not vanishing and the coefficient $a\left( x\right) 
$ is unknown, see \cite{Baud}.

While, in the case of CIPs, all above references are for the data resulting
from a single measurement event, in the paper \cite{KJIIP} the
data are generated by a point source running along an interval of a straight
line. A new version of the convexification method is constructed in \cite%
{KJIIP}. Detailed analytical and numerical studies of this version for the
case of the inverse problem of electrical impedance tomography were recently
conducted  in \cite{KEIT}.

The new version of the convexification for the 1-D CIP of this paper is
based on an adaptation of the idea of \cite{KJIIP} for the case of our 1-D
CIP. The governing PDE here is the 1-D Helmholtz equation, i.e. we work in
the frequency domain. One of important new elements of \cite{KJIIP} was the
construction of a new orthonormal basis in $L_{2}\left( 0,1\right) $ with
some special properties properties. In \cite{KJIIP} the functions of this
basis were dependent on the position of the point source, and the same was
in \cite{KEIT}.\ However, unlike \cite{KJIIP,KEIT}, we do not move the
source here.

Thus, a \emph{significantly new idea} of the current publication is that we
adapt that basis to the case when the wave number is running instead of the
point source in \cite{KJIIP}. In addition, unlike \cite{KlibKol1, KlibKol2},
the proof of the existence of the minimizer of our functional is different
here from the one of \cite{BakKlib}. Along with Carleman estimates, we use
here the apparatus of the convex analysis.

In addition to the new element mentioned in the preceding paragraph, we
eliminate here the following two restrictive conditions of \cite%
{KlibKol1,KlibKol2}:

\begin{enumerate}
\item Most importantly, unlike \cite{KlibKol1,KlibKol2}, we do not work here
with a nonlinear integro differential equation in which the integration is
carried out with respect to the wave number. This enables us to avoid the
use of the so-called \textquotedblleft tail function", which was used in 
\cite{KlibKol1,KlibKol2}. The tail function is the complement of a certain
truncated integral. Since the tail function is unknown, a certain
approximation for it was used in \cite{KlibKol1,KlibKol2}.

\item Unlike \cite{KlibKol1,KlibKol2}, we do not need to work here with
large values of wave numbers.
\end{enumerate}

CIPs have many applications in many different fields of science. Therefore,
there is a large variety of publications devoted to numerical solutions of
CIPs. We now refer to some of them as well as to the references cited
therein \cite{BK1, Ammari1, Ammari2, Ammari3, Chavent, Gonch1, Gonch2, Ito2012, Ito2013, Li2017}
. As to another globally convergent numerical method for some CIPs, we refer
to the technique of Kabanikhin and Shishlenin \cite{Kab1,Kab2,Kab3}. This
technique works with overdetermined data in the $n-$D case ($n=2,3$), unlike
the convexification.

In section 2 we state forward and inverse problems. In section 3 we
construct our weighted Tikhonov-like functional. In section 4 we formulate
our theorems. They are proved in section 5. In section 6 we present our
numerical results for both computationally simulated and experimental data.

\section{Statement of the inverse problem}

\label{sec:2}

Below $k>0$ is the wave number. Also, for any $z\in \mathbb{C}$ we denote $%
\overline{z}$ its complex conjugation.

Let $c_{0}>0$ be a positive number. Let $c\left( x\right) ,x\in \mathbb{R}$
be the function with the following properties: 
\begin{equation}
c\in C^{2}\left( \mathbb{R}\right), \quad c\left( x\right) \geq c_{0}, \quad \forall x\in 
\mathbb{R},  \label{2.1}
\end{equation}%
\begin{equation}
c\left( x\right) =1, \quad \forall x\notin \left( 0,1\right) ,  \label{2.2}
\end{equation}%
In our application $c\left( x\right) $ is the spatially distributed
dielectric constant of the medium. Let $x_{0}<0$ be the position of the
point source$.$ The forward problem is: 
\begin{equation}
u^{\prime \prime }+k^{2}c\left( x\right) u=-\delta \left( x-x_{0}\right), \quad x\in \mathbb{R},  \label{2.4}
\end{equation}%
\begin{equation}
\lim_{x\rightarrow \infty }\left( u^{\prime }+iku\right)
=0,\, \lim_{x\rightarrow -\infty }\left( u^{\prime }-iku\right) =0.  \label{2.6}
\end{equation}

\subsection{Statement of the inverse problem}

\label{sec:2.1}

Let $u_{0}\left( x,k\right) $ be the solution of the problem (\ref{2.4}), (%
\ref{2.6}) for the homogeneous case with $c\left( x\right) \equiv 1,$ 
\begin{equation}
u_{0}\left( x,k\right) =\frac{e^{-ik\left\vert x-x_{0}\right\vert }}{2ik}.
\label{2.60}
\end{equation}%
In this paper we consider the following coefficient inverse problem:

\textbf{Coefficient Inverse Problem (CIP)}.\emph{\ Let }$[\underline{k},%
\overline{k}]\subset \left( 0,\infty \right) $\emph{\ be an interval of wave
numbers }$k$\emph{. Reconstruct the function }$c\left( x\right) ,$\emph{\
assuming that the following function }$g_{0}\left( k\right) $\emph{\ is known%
} 
\begin{equation}
g_{0}\left( k\right) =\frac{u(0,k)}{u_{0}(0,k)},\quad k\in \lbrack \underline{k},%
\overline{k}].  \label{2.8}
\end{equation}

Uniqueness theorem for this CIP was proven in \cite{KlibLoc}. Denote%
\begin{equation}
w\left( x,k\right) =\frac{u\left( x,k\right) }{u_{0}\left( x,k\right) }.
\label{2.100}
\end{equation}%
By (\ref{2.8}) and (\ref{2.100}) 
\begin{equation}
w\left( 0,k\right) =g_{0}\left( k\right) ,\quad k\in \lbrack \underline{k},%
\overline{k}].  \label{2.101}
\end{equation}%
Besides, it was established in \cite{KlibLoc} that 
\begin{equation}
w^{\prime }\left( 0,k\right) =g_{1}\left( k\right) =2ik\left( g_{0}\left(
k\right) -1\right) ,\quad k\in \lbrack \underline{k},\overline{k}].  \label{2.160}
\end{equation}

\subsection{Some properties of the solution of the forward problem}

\label{sec:2.2}

Lemma 2.1 is formulated in \cite{KlibLoc} as Lemma 4.1. However, since the
proof was not provided in \cite{KlibLoc}, we prove Lemma 2.1 here.

\textbf{Lemma 2.1}. \emph{The function }$u\left( x,k\right) \neq 0$\emph{\
for all }$x>x_{0}$\emph{\ and for all }$k>0.$\emph{\ In particular, the
function }$g_{0}\left( k\right) \neq 0$\emph{\ for all }$k\in \lbrack 
\underline{k},\overline{k}].$

\textbf{Proof}. Since $x_{0}<0$ and by (\ref{2.2}) $c\left( x\right) =1$ for 
$x\geq 1,$ then (\ref{2.4}) implies that $u^{\prime \prime }+k^{2}u=0$ for $%
x\geq 1.$ Hence, 
\begin{equation}
u\left( x,k\right) =B_{1}\left( k\right) e^{ikx}+B_{2}\left( k\right)
e^{-ikx}, \quad  x\geq 1  \label{2.16}
\end{equation}%
with some $k-$dependent complex numbers $B_{1}\left( k\right) $ and $%
B_{2}\left( k\right) .$ Next, (\ref{2.6}) implies that 
\begin{equation}
B_{1}\left( k\right) =0.  \label{2.17}
\end{equation}

Assume that $u\left( x_{1},k\right) =0,$ where the point $x_{1}\in \left(
x_{0},1\right) .$ Multiply both sides of (\ref{2.4}) by $\overline{u}\left(
x,k\right) $ and integrate with respect to $x\in \left( x_{1},1\right) .$ We
obtain%
\begin{equation}
\left( u^{\prime }\overline{u}\right) \left( 1,k\right) -\displaystyle%
\int\limits_{x_{1}}^{1}\left\vert u^{\prime }\right\vert ^{2}dx+k^{2}%
\displaystyle\int\limits_{x_{1}}^{1}c\left( x\right) \left\vert u\right\vert
^{2}dx=0.  \label{2.18}
\end{equation}%
By (\ref{2.16}) and (\ref{2.17}) $u^{\prime }\left( 1,k\right) =-iku\left(
1,k\right) .$ Substituting this in (\ref{2.18}), we obtain%
\begin{equation}
-ik\left\vert u\left( 1,k\right) \right\vert ^{2}-\displaystyle%
\int\limits_{x_{1}}^{1}\left\vert u^{\prime }\right\vert ^{2}dx+k^{2}%
\displaystyle\int\limits_{x_{1}}^{1}c\left( x\right) \left\vert u\right\vert
^{2}dx=0.  \label{2.19}
\end{equation}%
Since the first term in the left hand side of (\ref{2.19}) is an imaginary
number, whereas two other terms are real numbers, then (\ref{2.19}) implies
that $u\left( 1,k\right) =0.$ Hence, (\ref{2.16}) and (\ref{2.17}) imply
that $u\left( x,k\right) =0$ for $x\geq 1.$ Hence, the uniqueness result for
the Cauchy problem for an ordinary ODE implies that 
\begin{equation}
u\left( x,k\right) =0\mbox{ for }x>x_{0}.  \label{2.190}
\end{equation}

Next, it follows from (\ref{2.2}) and (\ref{2.4}) that 
\begin{equation}
u\left( x,k\right) =B_{3}\left( k\right) e^{ikx},\quad  x<x_{0},  \label{2.191}
\end{equation}%
where the complex number $B_{3}\left( k\right) $ depends only on $k$. It is
well known (see, e.g. \cite{CK}) that the problem (\ref{2.4}), (\ref{2.6})
is equivalent with the 1-D Lippmann-Schwinger equation%
\begin{equation}
u\left( x,k\right) =\frac{e^{-ik\left\vert x-x_{0}\right\vert }}{2ik}+\frac{k%
}{2i}\int\limits_{0}^{1}e^{-ik\left\vert x-\xi \right\vert }\left( c\left(
\xi \right) -1\right) u\left( \xi ,k\right) d\xi .  \label{2.192}
\end{equation}%
Using (\ref{2.1}) and (\ref{2.192}), we obtain $u\in C\left( \mathbb{R}%
\right) \cap C^{3}\left( \mathbb{R}\diagdown \left\{ x_{0}\right\} \right) .$
Hence, (\ref{2.190}) implies that%
\begin{equation}
u\left( x_{0},k\right) =0.  \label{2.193}
\end{equation}%
Hence, setting in (\ref{2.191}) $x=x_{0},$ we obtain $B_{3}\left( k\right)
=0.$ The latter, (\ref{2.190}), (\ref{2.191}) and (\ref{2.193}) imply that $%
u\left( x,k\right) =0$ for all $x\in \mathbb{R}.$ This, however, is in the
contradiction with the fact that the right hand side of equation (\ref{2.4})
is not identical zero. $\square $

Lemma 2.1 helps us to define the function $\log w\left( x,k\right) .$ The
difficulty here is that the function $w\left( x,k\right) $ is a complex
valued one. Hence, the ambiguity of $\log w\left( x,k\right) $ might be up
to $2n\pi i,$ where $n$ is an integer. When defining $\log w\left(
x,k\right) ,$ we follow \cite{KlibLoc}. Consider the function $\phi \left(
x\right) ,$ 
\[
\phi \left( x\right) =-\frac{c^{\prime \prime }\left( x\right) }{c^{2}\left(
x\right) }+\frac{7}{16}\frac{\left( c^{\prime }\left( x\right) \right) ^{2}}{%
c^{3}\left( x\right) }. 
\]%
It was proven in theorem 2.2 of \cite{KlibLoc} that if 
\begin{equation}
\phi \left( x\right) \leq 0,  \label{2.20}
\end{equation}%
then for $x>x_{0}$ the asymptotic behavior of the function $u\left(
x,k\right) $ at $k\rightarrow \infty $ is 
\begin{equation}
u\left( x,k\right) =\frac{1}{2ikc^{1/4}\left( x\right) }\exp \left[ -ik%
\displaystyle\int\limits_{x_{0}}^{x}\sqrt{c\left( \xi \right) }d\xi \right]
\left( 1+O\left( \frac{1}{k}\right) \right) .  \label{2.21}
\end{equation}%
Assuming that this asymptotic behavior holds regardless on the inequality (%
\ref{2.20}), we now remind the definition of \cite{KlibLoc} of the function $%
\log w\left( x,k\right) .$ First, it follows from (\ref{2.60}), (\ref{2.100}%
) and (\ref{2.21}) that we can define $\log w\left( x,k\right) $ for
sufficiently large values of $k$ and for $x\in \left( 0,1\right) $ as 
\begin{equation}
\log w\left( x,k\right) =-\frac{1}{4}\ln c\left( x\right) -ik\left(
\int\limits_{x_{0}}^{x}\sqrt{c\left( \xi \right) }d\xi -x+x_{0}\right)
\label{2.22}
\end{equation}%
\[
-\ln \left( 2k\right) -i\frac{\pi }{2}+\log \left( 1+\widehat{w}\left(
x,k\right) \right) , 
\]%
where 
\begin{equation}
\widehat{w}\left( x,k\right) =O\left( \frac{1}{k}\right) ,\quad k\rightarrow
\infty .  \label{2.23}
\end{equation}%
Hence, we set for sufficiently large $k$: 
\begin{equation}
\log \left( 1+\widehat{w}\left( x,k\right) \right) =\displaystyle%
\sum\limits_{n=1}^{\infty }\left( -1\right) ^{n-1}\frac{\left( \widehat{w}%
\left( x,k\right) \right) ^{n}}{n},  \label{2.24}
\end{equation}%
which eliminates the above mentioned ambiguity.

Let $K>0$ be such a number that (\ref{2.24}) is valid for $k\geq K.$ Note
that $K$ might depend on $x$, $K=K\left( x\right) $. For $k\in \left(
0,K\right) $ consider the function $\psi (x,k),$ 
\begin{equation}
\psi (x,k)=-\int\limits_{k}^{\overline{k}}\frac{\partial _{k}w(x,\kappa )}{%
w(x,\kappa )}d\kappa +\log w(x,K),  \label{2.25}
\end{equation}%
where $\log w(x,K)$ is defined via (\ref{2.22}) and (\ref{2.24}). By (\ref%
{2.100}) and Lemma 2.1 $w(x,\kappa )\neq 0,\forall x\in \left[ 0,1\right]
,\forall \kappa >0.$ Differentiate both sides of (\ref{2.25}) with respect
to $k$. We obtain 
\begin{equation}
\partial _{k}w(x,k)-w(x,k)\partial _{k}\psi (x,k)=0.  \label{2.26}
\end{equation}%
Multiplying both sides of (\ref{2.26}) by $\exp (-\psi (x,k))$, we obtain $%
\partial _{k}\left( e^{-\psi (x,k)}w(x,k)\right) =0.$ Hence, there exists a
function $C=C\left( x\right) $ independent on $k$ such that 
\begin{equation}
w(x,k)=C\left( x\right) e^{\psi (x,k)}.  \label{2.27}
\end{equation}%
Set in (\ref{2.27}) $k=K.$ By (\ref{2.25}) $\psi (x,K)=\log w(x,K)$. Hence,
we obtain from (\ref{2.27}) 
\begin{equation}
C=C\left( x\right) =1, \quad x\in \left[ 0,1\right] .  \label{2.28}
\end{equation}%
Hence, it follows from (\ref{2.27}) and (\ref{2.28}) that it is appropriate
to define $\log w(x,k)$ as $\log w(x,k)=\psi (x,k).$ And ambiguity does
occur this way.

\section{The Tikhonov-like Functional With the CWF}

\label{sec:3}

In this section we construct the above mentioned weighted cost functional
with the CWF in it.

\subsection{A special orthonormal basis in $L_{2}\left( \protect\underline{k}%
,\overline{k}\right) $}

\label{sec:3.1}

This basis was first constructed by Klibanov in \cite{KJIIP}. In this
reference, functions of that basis are dependent on the position of the
point source, and this source runs along an interval of a straight line. We
now briefly repeat that construction for the case when these functions
depend on the wave number $k$. The latter is a new idea.

Let $\left( ,\right) $ denotes the scalar product in $L_{2}\left( \underline{%
k},\overline{k}\right) .$ We need to construct an orthonormal basis $\left\{
\psi _{n}\left( k\right) \right\} _{n=0}^{\infty }$ of real valued function
in the space $L_{2}\left( \underline{k},\overline{k}\right) $ such that the
following two conditions are met:

\begin{enumerate}
\item $\psi _{n}\in C^{1}\left[ \underline{k},\overline{k}\right] ,$ $%
\forall n=0,1,...$

\item Let $a_{mn}=\left( \psi _{n}^{\prime },\psi _{m}\right) .$ Then the
matrix $M_{N}=\left( a_{mn}\right) _{m,n=0}^{N-1}$ is invertible for any $%
N=1,2,...$
\end{enumerate}

Neither classical orthonormal polynomials nor the basis of trigonometric
functions $\left\{ \exp \left[ in\pi \left( k-\underline{k}\right) /\left( 
\overline{k}-\underline{k}\right) \right] \right\} _{n=0}^{\infty }$ do not
satisfy the second condition. This is because in any of these two cases the
first column of the matrix $M_{N}$ would be equal to zero.

Consider the set of functions $\left\{ k^{n}e^{k}\right\} _{n=0}^{\infty }.$
This set is complete in $L_{2}\left( 0,1\right) .$ We orthonormalize it
using the classical Gram-Schmidt orthonormalization procedure. We start from 
$n=0$, then take $n=1$, etc. Then we obtain the orthonormal basis $\left\{
\varphi _{n}\left( k\right) \right\} _{n=0}^{\infty }$ in $L_{2}\left(
0,1\right) .$ Each function $\varphi _{n}\left( k\right) $ has the form%
\[
\varphi _{n}\left( k\right) =p_{n}\left( k\right) e^{k}, 
\]%
where $p_{n}\left( k\right) $ is the polynomial of the degree $n$. Next, we
set 
\[
\psi _{n}\left( k\right) =\frac{1}{\sqrt{\overline{k}-\underline{k}}}\varphi
_{n}\left( \frac{k-\underline{k}}{\overline{k}-\underline{k}}\right) . 
\]%
Thus, the set $\left\{ \psi _{n}\left( k\right) \right\} _{n=0}^{\infty }$
is an orthonormal basis in the space $L_{2}\left( \underline{k},\overline{k}%
\right) .$ Lemma 3.1 ensures that the above property number 2 holds for
functions $\psi _{n}\left( k\right) $.

\textbf{Lemma 3.1} \cite{KJIIP}\textbf{.} \emph{Denote }$a_{mn}=\left( \psi
_{n}^{\prime },\psi _{m}\right) .$ \emph{Then} 
\begin{equation}
a_{mn}=\left\{ 
\begin{array}{c}
\left( \overline{k}-\underline{k}\right) ^{-1}\mbox{ if }n=m, \\ 
0\mbox{ if }n<m.%
\end{array}%
\right.  \label{3.2}
\end{equation}%
\emph{For an integer }$N\geq 1$\emph{\ consider the }$N\times N$\emph{\
matrix }$M_{N}=\left( a_{mn}\right) _{\left( m,n\right) =\left( 0,0\right)
}^{\left( N-1,N-1\right) }.$\emph{\ Then (\ref{3.2}) implies that }$M_{N}$%
\emph{\ is an upper diagonal matrix with} $\det \left( M_{N}\right) =\left( 
\overline{k}-\underline{k}\right) ^{-N}\neq 0.$ \emph{Thus, the inverse
matrix }$M_{N}^{-1}$\emph{\ exists. }

\subsection{A system of coupled quasilinear ordinary differential equations}

\label{sec:3.2}

Let $\log w\left( x,k\right) $ be the function constructed in section 2.2.
Consider the function%
\begin{equation}
v\left( x,k\right) =\frac{\log w\left( x,k\right) }{k^{2}},\,\quad x\in \left[ 0,1%
\right] , \, k\in \left[ \underline{k},\overline{k}\right] .  \label{eq:v}
\end{equation}%
Substitution in equation (\ref{2.4}) leads to%
\begin{equation}
v^{\prime \prime }+k^{2}\left( v^{\prime }\right) ^{2}-2ikv^{\prime
}=1-c\left( x\right) =-\beta (x), \quad x\in \left[ 0,1\right] , \, k\in \left[ 
\underline{k},\overline{k}\right] .  \label{3.3}
\end{equation}%
Also, by (\ref{2.100})-(\ref{2.160})%
\begin{equation}
v\left( 0,k\right) =q_{0}\left( k\right) , \quad  v^{\prime }\left( 0,k\right)
=q_{1}\left( k\right) ,  \label{3.4}
\end{equation}%
where 
\begin{equation}
q_{0}\left( k\right) =\frac{\log g_{0}\left( k\right) }{k^{2}},\quad q_{1}\left(
k\right) =\frac{g_{1}\left( k\right) }{g_{0}\left( k\right) k^{2}}.
\label{eq:v_bcs}
\end{equation}%
In addition, (\ref{2.16}) and (\ref{2.17}) imply that 
\begin{equation}
v^{\prime }\left( 1,k\right) =0.  \label{3.5}
\end{equation}%
Differentiate both sides of equation (\ref{3.3}) with respect to $k$ and use
the fact that $\partial _{k}\left( 1-c\left( x\right) \right) =0.$ We obtain%
\begin{equation}
v_{k}^{\prime \prime }+2k^{2}v_{k}^{\prime }v^{\prime }+2k\left( v^{\prime
}\right) ^{2}-2ikv_{k}^{\prime }-2iv^{\prime }=0.  \label{3.6}
\end{equation}

We now assume that the function $v\left( x,k\right) $ can be represented via
a truncated Fourier series with respect to the orthonormal basis $\left\{
\psi _{n}\left( k\right) \right\} _{n=0}^{\infty }$. In fact, this is our
main approximation. Thus, for an integer $N\geq 1,$ we assume that 
\begin{equation}
v\left( x,k\right) =\displaystyle\sum\limits_{n=0}^{N-1}y_{n}\left( x\right)
\psi _{n}\left( k\right)  \label{3.7}
\end{equation}%
for $x\in \left[ 0,1\right] ,k\in \left[ \underline{k},\overline{k}\right] .$
Here coefficients $y_{n}\left( x\right) $ of the expansion (\ref{3.7}) are
unknown and should be determined as the main part of the solution of our
CIP. In particular, (\ref{3.3}) and (\ref{3.7})\ imply that%
\begin{equation}
c\left( x\right) =1-\sum\limits_{n=0}^{N-1}y_{n}^{\prime \prime }\left(
x\right) \psi _{n}\left( k\right) -k^{2}\left( \displaystyle%
\sum\limits_{n=0}^{N-1}y_{n}^{\prime }\left( x\right) \psi _{n}\left(
k\right) \right) ^{2}   \label{3.8}
\end{equation}%
\[
+2ik\sum\limits_{n=0}^{N-1}y_{n}^{\prime }\left( x\right) \psi _{n}\left(
k\right) . 
\]%
Substituting (\ref{3.7}) in (\ref{3.6}), we obtain%
\[
\displaystyle\sum\limits_{n=0}^{N-1}y_{n}^{\prime \prime }\left( x\right)
\psi _{n}^{\prime }\left( k\right) +2k^{2}\displaystyle\sum%
\limits_{n,m=0}^{N-1}y_{n}^{\prime }\left( x\right) y_{m}^{\prime }\left(
x\right) \psi _{n}^{\prime }\left( k\right) \psi _{m}^{\prime }\left(
k\right) 
\]%
\begin{equation}
+2k\left[ \displaystyle\sum\limits_{n=0}^{N-1}y_{n}^{\prime }\left( x\right)
\psi _{n}\left( k\right) \right] ^{2}-2ik\displaystyle\sum%
\limits_{n=0}^{N-1}y_{n}^{\prime }\left( x\right) \psi _{n}^{\prime }\left(
k\right)  \label{3.9}
\end{equation}%
\[
-2i\displaystyle\sum\limits_{n=0}^{N-1}y_{n}^{\prime }\left( x\right) \psi
_{n}\left( k\right) =0. 
\]%
Introduce the $N-$D vector function $y\left( x\right) ,$ 
\begin{equation}
y\left( x\right) =\left( y_{0},...,y_{N-1}\right) ^{T}\left( x\right) .
\label{3.90}
\end{equation}%
Let $s\in \left[ 0,N-1\right] $ be an integer. Multiply both sides of (\ref%
{3.9}) by $\psi _{s}\left( k\right) $ and integrate with respect to $k\in
\left( \underline{k},\overline{k}\right) .$ We obtain%
\begin{equation}
M_{N}y^{\prime \prime }+\widetilde{F}\left( y^{\prime }\right) =0,\quad x\in %
\left[ 0,1\right] .  \label{3.10}
\end{equation}%
Here the vector function $\widetilde{F}\left( \cdot \right) $ is quadratic
with respect to the functions $y_{n}^{\prime }\left( x\right) .$ Since by
Lemma 3.1 the matrix $M_{N}$ is invertible, we multiply both sides of (\ref%
{3.10}) by $M_{N}^{-1}$ and obtain%
\begin{equation}
y^{\prime \prime }+F\left( y^{\prime }\right) =0, \quad x\in \left[ 0,1\right] ,
\label{3.11}
\end{equation}%
where the vector function $F\left( y^{\prime }\right) =M_{N}^{-1}\widetilde{F%
}\left( y^{\prime }\right) .$ Therefore, the vector function $F\left(
y^{\prime }\right) $ is quadratic with respect to the functions $%
y_{n}^{\prime }\left( x\right) .$ In addition, (\ref{3.4}) and (\ref{3.5})
imply that $y\left( x\right) $ satisfies the following boundary conditions:%
\begin{equation}
y\left( 0\right) =f_{0}, \quad y^{\prime }\left( 0\right) =f_{1}, \quad y^{\prime }\left(
1\right) =0.  \label{3.12}
\end{equation}%
Note that due to the boundary condition $y^{\prime }\left( 1\right) =0,$ (%
\ref{3.11}), (\ref{3.12}) is not the regular Cauchy problem for the system (%
\ref{3.11}) of coupled quasilinear Ordinary Differential Equations. Thus,
our effort below is focused on the numerical solution of the problem (\ref%
{3.11}), (\ref{3.12}). Indeed, if we would solve it, then, using (\ref{3.8})
and (\ref{3.90}), we would find the target coefficient $c\left( x\right) .$
We solve this problem via minimizing our weighted Tikhonov-like functional
with the CWF. First, we formulate the Carleman estimate for the
operator $d^{2}/dx^{2}.$

\textbf{Lemma 3.2} (Carleman estimate) \cite{KlibKol1}. \emph{For any
complex valued function }$u\in H^{2}\left( 0,1\right) $\emph{\ with }$%
u(0)=u^{\prime }(0)=0$\emph{\ and for any parameter }$\lambda \geq 1$\emph{\
the following Carleman estimate holds }%
\[
\int\limits_{0}^{1}\left\vert u^{\prime \prime }\right\vert ^{2}e^{-2\lambda
x}dx\geq C\int\limits_{0}^{1}|u^{\prime \prime }|^{2}e^{-2\lambda x}dx 
\]%
\[
+C\left[ \lambda \int\limits_{0}^{1}|u^{\prime }|^{2}e^{-2\lambda
x}dx+\lambda ^{3}\int\limits_{0}^{1}|u|^{2}e^{-2\lambda x}dx\right] , 
\]%
\emph{where\ the constant }$C>0$ \emph{is} \emph{independent on }$u$\emph{\
and }$\lambda .$

\textbf{Remark 3.1}. \emph{Using Lemma 3.1 as well as arguments, which are
completely similar with those of the proof of Theorem 3.1 of \cite{KJIIP},
one can prove that there exists at most one solution }$V\in H^{2}\left(
0,1\right) $\emph{\ of the problem (\ref{3.11}), (\ref{3.12}).}

\subsection{ Weighted Tikhonov-like functional}

\label{sec:3.3}

Introduce the cut-off function $\chi \left( x\right) ,$%
\begin{equation}
\chi \in C^{2}\left[ 0,1\right] , \quad \chi \left( x\right) =\left\{ 
\begin{array}{c}
1,x\in \left[ 0,1/2\right] , \\ 
0,x\in \left[ 3/4,1\right] , \\ 
\in \left( 0,1\right) ,x\in \left( 1/2,3/4\right) .%
\end{array}%
\right.  \label{1}
\end{equation}%
Consider the $N-$D vector function $f\left( x\right) $ defined as%
\begin{equation}
f\left( x\right) =\left[ f_{0}+xf_{1}\right] \chi \left( x\right) .
\label{2}
\end{equation}%
Then%
\[
f\in C^{2}\left[ 0,1\right] , \quad f\left( 0\right) =f_{0}, \quad f^{\prime }\left(
0\right) =f_{1}, \quad f^{\prime }(1)=0. 
\]%
Hence%
\begin{equation}
\left( y-f\right) \left( 0\right) =0,\quad  \left( y-f\right) ^{\prime }\left(
0\right) =0, \quad \left( y-f\right) ^{\prime }\left( 1\right) =0.  \label{eq:Vf}
\end{equation}

Let $R>0$ be an arbitrary number. Consider the sets $B\left(
R,f_{0},f_{1}\right) $ and $B_{0}\left( R\right) $ of $N-$D vector functions 
$W\left( x\right) $ and $p\left( x\right) $ defined as:%
\begin{equation}
B\left( R,f_{0},f_{1}\right) =\left\{ 
\begin{array}{c}
W\in H^{2}\left( 0,1\right) : \\ 
W\left( 0\right) =f_{0},W^{\prime }\left( 0\right) =f_{1},W^{\prime }\left(
1\right) =0, \\ 
\left\Vert W-f\right\Vert _{H^{2}\left( 0,1\right) }<R%
\end{array}%
\right\} ,  \label{3.13}
\end{equation}%
\begin{equation}
B_{0}\left( R\right) =\left\{ p\in H^{2}\left( 0,1\right) :p\left( 0\right)
=p^{\prime }\left( 0\right) =p^{\prime }\left( 1\right) =0,\left\Vert
p\right\Vert _{H^{2}\left( 0,1\right) }<R\right\} .  \label{3.132}
\end{equation}%
Hence, all vector functions belonging to the set $B\left(
R,f_{0},f_{1}\right) $ satisfy boundary conditions (\ref{3.12}). Note that
both sets $B\left( R,f_{0},f_{1}\right) $ and $B_{0}\left( R\right) $ are
convex. The proof of Proposition 3.1 follows immediately from (\ref{3.13})
and (\ref{3.132}).

\textbf{Proposition 3.1.} \emph{Let }$f\left( x\right) $\emph{\ be the
vector function defined in (\ref{2}).} \emph{Then} \emph{for any vector
function }$W\in B\left( R,f_{0},f_{1}\right) $\emph{\ the vector function }$%
\left( W-f\right) \in B_{0}\left( R\right) .$\emph{\ And vice versa: for any
vector function }$p\in B_{0}\left( R\right) $\emph{\ the vector function }$%
\left( p+f\right) \in B\left( R,f_{0},f_{1}\right) .$

In accordance with the Tikhonov concept for ill-posed problems \cite{Bak},
we assume that there exists the exact solution $y^{\ast }$ of the problem (%
\ref{3.11}), (\ref{3.12}) with the exact (i.e. noiseless) data 
\begin{equation}
y^{\ast }\left( 0\right) =f_{0}^{\ast },\quad y^{\ast \prime }\left( 0\right)
=f_{1}^{\ast }  \label{4.3}
\end{equation}%
in (\ref{3.12}) and $y^{\ast }\in B\left( R,f_{0}^{\ast },f_{1}^{\ast
}\right) .$ Let a sufficiently small number $\delta \in \left( 0,1\right) $
represents the level of the noise in the data, i.e.%
\begin{equation}
\left\vert f_{0}-f_{0}^{\ast }\right\vert <\delta ,\quad \left\vert
f_{1}-f_{1}^{\ast }\right\vert <\delta .  \label{4.4}
\end{equation}%
Similarly with (\ref{2}) introduce the function $f^{\ast }\left( x\right) ,$%
\begin{equation}
f^{\ast }\left( x\right) =\left[ f_{0}^{\ast }+xf_{1}^{\ast }\right] \chi
\left( x\right) .  \label{3.130}
\end{equation}%
It follows from (\ref{1}), (\ref{2}), (\ref{4.4}) and (\ref{3.130}) that%
\begin{equation}
\left\Vert f-f^{\ast }\right\Vert _{C^{2}\left[ 0,1\right] }\leq B\delta .
\label{3.131}
\end{equation}%
Here and below $B=B\left( \chi \right) >0$ denotes different positive
numbers depending only on the function $\chi \left( x\right) .$ Since $%
\delta $ is sufficiently small, we indicate below dependencies of some
constants on $f^{\ast }$ rather than on $f$.

Let $\alpha \in \left( 0,1\right) $ be the regularization parameter. Our
Tikhonov-like weighted functional is%
\begin{equation}
J_{\lambda ,\alpha }\left( y\right) =e^{2\lambda
}\int\limits_{0}^{1}\left\vert y^{\prime \prime }+F\left( y^{\prime }\right)
\right\vert ^{2}e^{-2\lambda x}dx+\alpha \left\Vert y\right\Vert
_{H^{2}\left( 0,1\right) }^{2}.  \label{3.14}
\end{equation}%
Here the multiplier $e^{2\lambda }$ is introduced to balance two terms in
the right hand side of (\ref{3.14}). We consider the following minimization
problem:

\textbf{Minimization Problem}. \emph{Minimize the functional }$J_{\lambda
,\alpha }\left( y\right) $\emph{\ on the set }$y\in \overline{B\left(
R,f_{0},f_{1}\right) }.$

\textbf{Remark 3.1}. \emph{In principle, estimate (\ref{4.1}) of Theorem 4.1
tells one that it is not necessary to incorporate the regularization term }$%
\alpha \left\Vert y\right\Vert _{H^{2}\left( 0,1\right) }^{2}$\emph{\ in the
functional }$J_{\lambda ,\alpha }\left( y\right) .$\emph{\ Nevertheless we
have observed in our computations that the presence of this term improves
numerical results. This is why we introduce it here. We cannot yet explain
the reason of this improvement. }

\section{Theorems}

\label{sec:4}

In this section we formulate theorems about the functional $J_{\lambda
,\alpha }\left( y\right) .$ We prove them in section 5. Theorem 4.1 is the
central result of this paper, also see Remark 3.1 about $\alpha $. Below $%
C_{1}=C_{1}\left( R,F,N\right) >0$ and $C_{2}=C_{2}\left( R,F,N,\chi
,f^{\ast }\right) >0$ denote different numbers depending only on listed
parameters.

\textbf{Theorem 4.1.} \emph{The functional }$J_{\lambda ,\alpha }\left(
y\right) $\emph{\ has the Frech\'{e}t derivative }$J_{\lambda ,\alpha
}^{\prime }\left( y\right) $\emph{\ at each point }$y\in B\left(
2R,f_{0},f_{1}\right) .$\emph{\ Also, there exists a number }$\lambda
_{1}=\lambda _{1}\left( R,F,N\right) >1$\emph{\ depending only on listed
parameters such that for all }$\lambda \geq \lambda _{1}$\emph{\ the
functional }$J_{\lambda ,\alpha }\left( y\right) $\emph{\ is strictly convex
on the set }$\overline{B\left( R,f_{0},f_{1}\right) },$\emph{\ i.e. for all }%
$y_{\left( 1\right) },y_{\left( 2\right) }\in \overline{B\left(
R,f_{0},f_{1}\right) }$%
\begin{equation}
J_{\lambda ,\alpha }\left( y_{\left( 2\right) }\right) -J_{\lambda ,\alpha
}\left( y_{\left( 1\right) }\right) -J_{\lambda ,\alpha }^{\prime }\left(
y_{\left( 1\right) }\right) \left( y_{\left( 2\right) }-y_{\left( 1\right)
}\right) \geq C_{1}\left\Vert y_{\left( 2\right) }-y_{\left( 1\right)
}\right\Vert _{H^{2}\left( 0,1\right) }^{2}.  \label{4.1}
\end{equation}

\textbf{Corollary 4.1. }\emph{Consider the functional }%
\begin{equation}
\Phi _{\lambda ,\alpha }\left( p\right) =J_{\lambda ,\alpha }\left(
p+f\right) , \quad \forall p\in \overline{B_{0}\left( R\right) }.  \label{5.22}
\end{equation}%
\emph{Then a direct analog of Theorem 4.1 is valid for this functional. More
precisely, the functional }$\Phi _{\lambda ,\alpha }\left( P\right) $\emph{\
has the Frech\'{e}t derivative }$\Phi _{\lambda ,\alpha }^{\prime }\left(
p\right) $\emph{\ at each point }$p\in B_{0}\left( 2R\right) .$\emph{\ Let }$%
\lambda _{1}=\lambda _{1}\left( R,F,N\right) >1$\emph{\ be the number of
Theorem 4.1. Then there exists a number }$\lambda _{2}=\lambda _{2}\left(
R,F,N,\chi ,f^{\ast }\right) \geq \lambda _{1}$\emph{\ depending only on
listed parameters such that for all }$\lambda \geq \lambda _{2}$\emph{\ the
functional }$\Phi _{\lambda ,\alpha }\left( p\right) $\emph{\ is strictly
convex on the set }$\overline{B_{0}\left( R\right) },$\emph{\ i.e. for all }$%
p_{1},p_{2}\in \overline{B_{0}\left( R\right) }$%
\[
\Phi _{\lambda ,\alpha }\left( p_{2}\right) -\Phi _{\lambda ,\alpha }\left(
p_{1}\right) -\Phi _{\lambda ,\alpha }^{\prime }\left( p_{1}\right) \left(
p_{2}-p_{1}\right) \geq C_{2}\left\Vert p_{2}-p_{1}\right\Vert _{H^{2}\left(
0,1\right) }^{2}. 
\]

\textbf{Theorem 4.2}. \emph{The Frech\'{e}t derivatives }$J_{\lambda ,\alpha
}^{\prime }\left( y\right) $\emph{\ and }$\Phi _{\lambda ,\alpha }^{\prime
}\left( P\right) $ \emph{of both functionals }$J_{\lambda ,\alpha }\left(
y\right) $\emph{\ and }$\Phi _{\lambda ,\alpha }\left( P\right) $ \emph{are
Lipschitz continuous on }$B\left( 2R,f_{0},f_{1}\right) $ \emph{and }$%
B_{0}\left( 2R\right) $\emph{\ respectively.\ In other words, there exists a
constant }$D=D\left( R,F,\lambda ,\alpha \right) >0$\emph{\ depending only
on listed parameters such that for all }$y_{\left( 1\right) },y_{\left(
2\right) }\in B\left( 2R,f_{0},f_{1}\right) $\emph{\ }%
\[
\left\Vert J_{\lambda ,\alpha }^{\prime }\left( y_{\left( 1\right) }\right)
-J_{\lambda ,\alpha }^{\prime }\left( y_{\left( 2\right) }\right)
\right\Vert _{H^{2}\left( 0,1\right) }\leq D\left\Vert y_{\left( 2\right)
}-y_{\left( 1\right) }\right\Vert _{H^{2}\left( 0,1\right) } 
\]%
\emph{and similarly for} $\Phi _{\lambda ,\alpha }^{\prime }\left( P\right)
. $

\textbf{Theorem 4.3}. \emph{Let }$\lambda _{1}=\lambda _{1}\left(
R,F,N\right) >1$\emph{\ and }$\lambda _{2}=\lambda _{2}\left( R,F,N,\chi
,f^{\ast }\right) \geq \lambda _{1}$ \emph{be the numbers of Theorem 4.1 and
Corollary 4.1 respectively. Let }$f\left( x\right) $\emph{\ be the function
defined in (\ref{1}), (\ref{2}). Then for any }$\lambda \geq \lambda _{2}$%
\emph{\ and for any} $\alpha \in \left( 0,1\right) $ \emph{there exists a
unique minimizer }$y_{\min ,\lambda ,\alpha }$\emph{\ of the functional }$%
J_{\lambda ,\alpha }\left( y\right) $\emph{\ on the set }$\overline{B\left(
R,f_{0},f_{1}\right) }.$\emph{\ In addition, for these values of }$\lambda $ 
\emph{and }$\alpha $\emph{\ there exists unique minimizer }$p_{\min ,\lambda
,\alpha }$\emph{\ of the functional }$\Phi _{\lambda ,\alpha }\left(
p\right) $\emph{\ on the set }$\overline{B_{0}\left( R\right) }$\emph{.
Furthermore, }$p_{\min ,\lambda ,\alpha }=$\emph{\ }$y_{\min ,\lambda
,\alpha }-f$ \emph{and}%
\begin{equation}
\Phi _{\lambda ,\alpha }^{\prime }\left( p_{\min ,\lambda ,\alpha }\right)
\left( p_{\min ,\lambda ,\alpha }-p\right) \leq 0,\quad \forall p\in \overline{%
B_{0}\left( R\right) },  \label{4.20}
\end{equation}%
\emph{\ }%
\begin{equation}
J_{\lambda ,\alpha }^{\prime }\left( y_{\min ,\lambda ,\alpha }\right)
\left( y_{\min ,\lambda ,\alpha }-y\right) \leq 0, \quad \forall y\in \overline{%
B\left( R,f_{0},f_{1}\right) }.  \label{4.2}
\end{equation}

\textbf{Theorem 4.4} (accuracy estimate). \emph{Assume that the exact
solution }$y^{\ast }$ \emph{of the problem (\ref{3.11}), (\ref{3.12}) exists
and }$y^{\ast }\in B\left( R,f_{0}^{\ast },f_{1}^{\ast }\right) .$\emph{\
Also, assume that (\ref{4.3}) and (\ref{4.4}) hold. Denote }$p^{\ast
}=y^{\ast }-f^{\ast }\in B_{0}\left( R\right) .$ \emph{In addition, assume
that there exists the exact solution }$c^{\ast }\left( x\right) $\emph{\ of
our CIP and this function satisfies conditions (\ref{2.1}), (\ref{2.2}).
Suppose that the function }$c^{\ast }\left( x\right) $\emph{\ can be found
from the vector function }$y^{\ast }\left( x\right) $\emph{\ via formula (%
\ref{3.8}) in which functions }$y_{n}$\emph{\ are replaced with components }$%
y_{n}^{\ast }$\emph{\ of the vector function }$y^{\ast }\left( x\right) .$%
\emph{\ Let }$\lambda _{1}=\lambda _{1}\left( R,F,N\right) >1$\emph{\ and }$%
\lambda _{2}=\lambda _{2}\left( R,F,N,\chi ,f^{\ast }\right) \geq \lambda
_{1}$ \emph{be the numbers of Theorem 4.1 and Corollary 4.1 respectively.
Consider the number }$\delta _{0}$ \emph{such that }$\delta _{0}\in \left(
0,e^{-4\lambda _{2}}\right) .$\emph{\ For any }$\delta \in \left( 0,\delta
_{0}\right) $\emph{\ we set} $\lambda =\lambda \left( \delta \right) =\ln
\left( \delta ^{-1/4}\right) >\lambda _{2}$ \emph{and} $\alpha =\alpha
\left( \delta \right) =\sqrt{\delta }.$ \emph{Let }$y_{\min ,\lambda ,\alpha
}\in B\left( R,f_{0},f_{1}\right) $\emph{\ and }$p_{\min ,\lambda ,\alpha
}=y_{\min ,\lambda ,\alpha }-f$\emph{\ }$\in B_{0}\left( R\right) \ $\emph{%
be the unique minimizers of the functionals }$J_{\lambda ,\alpha }\left(
V\right) $\emph{\ and }$\Phi _{\lambda ,\alpha }\left( p\right) $\emph{\ on
sets }$B\left( R,f_{0},f_{1}\right) $\emph{\ and }$B_{0}\left( R\right) $%
\emph{\ respectively (Theorem 4.3). Then the following accuracy estimates
hold: }%
\begin{equation}
\left\Vert p_{\min ,\lambda ,\alpha }-p^{\ast }\right\Vert _{H^{2}\left(
0,1\right) }\leq C_{2}\delta ^{1/4},  \label{4.5}
\end{equation}%
\begin{equation}
\left\Vert y_{\min ,\lambda \left( \delta \right) ,\alpha \left( \delta
\right) }-y^{\ast }\right\Vert _{H^{2}\left( 0,1\right) }\leq C_{2}\delta
^{1/4},  \label{4.50}
\end{equation}%
\begin{equation}
\left\Vert c_{\min ,\lambda \left( \delta \right) ,\alpha \left( \delta
\right) }-c^{\ast }\right\Vert _{L_{2}\left( 0,1\right) }\leq C_{2}\delta
^{1/4},  \label{4.6}
\end{equation}%
\emph{where the function }$c_{\min ,\lambda \left( \delta \right) ,\alpha
\left( \delta \right) }$\emph{\ is found from components }$y_{n,\min
,\lambda \left( \delta \right) ,\alpha \left( \delta \right) }$\emph{\ of
the vector function }$y_{\min ,\lambda \left( \delta \right) ,\alpha \left(
\delta \right) }$\emph{\ via formula (\ref{3.8}).}

Define the subspace $H_{0}^{2}\left( 0,1\right) $ of the space $H^{2}\left(
0,1\right) $ as%
\[
H_{0}^{2}\left( 0,1\right) =\left\{ w\in H^{2}\left( 0,1\right) :w\left(
0\right) =w^{\prime }\left( 0\right) =w^{\prime }\left( 1\right) =0\right\}
. 
\]%
By (\ref{3.132}) $B_{0}\left( R\right) \subset H_{0}^{2}\left( 0,1\right) .$
Let $Q_{\overline{B_{0}}}:H_{0}^{2}\left( 0,1\right) \rightarrow \overline{%
B_{0}\left( R\right) }$ be the projection operator of the space $%
H_{0}^{2}\left( 0,1\right) $ on the closed ball $\overline{B_{0}\left(
R\right) }.$ Consider now the gradient projection method of the minimization
of the functional (\ref{5.22}) on the set $\overline{B_{0}\left( R\right) }.$%
\ Let $p_{0}\in B_{0}\left( R\right) $ be an arbitrary point and $\gamma >0$
be a number. We consider the following sequence:%
\begin{equation}
p_{n}=Q_{\overline{B_{0}}}\left( p_{n-1}-\gamma \Phi _{\lambda ,\alpha
}^{\prime }\left( p_{n-1}\right) \right) ,\quad n=1,2,...  \label{4.7}
\end{equation}

\textbf{Theorem 4.5} (global convergence of the gradient method). \emph{Let
conditions of Theorem 4.4 about exact solutions }$y^{\ast }$\emph{\ and }$%
c^{\ast }$\emph{\ hold.\ Let the numbers }$\lambda _{1},\lambda _{2},\delta
,\lambda \left( \delta \right) ,\alpha \left( \delta \right) $\emph{\ be the
same as in Theorem 4.4. Let an arbitrary point }$p_{0}\in B_{0}\left(
R\right) $\emph{\ be} \emph{the starting point of the gradient projection
method (\ref{4.7})}$.$ \emph{Consider the sequence (\ref{4.7})} \emph{and
denote }$y^{n}=p_{n}+f.$ \emph{Then there exists a number }$\gamma
_{0}=\gamma _{0}\left( R,F,N,\chi ,f^{\ast },\delta \right) \in \left(
0,1\right) $\emph{\ depending only on listed parameters such that for any }$%
\gamma \in \left( 0,\gamma _{0}\right) $\emph{\ there exists a number }$%
q=q\left( \gamma \right) \in \left( 0,1\right) $ such that\emph{\ the
following convergence estimates are valid:}%
\begin{equation}
\left\Vert p_{n}-p_{\min ,\lambda \left( \delta \right) ,\alpha \left(
\delta \right) }\right\Vert _{H^{2}\left( 0,1\right) }\leq q^{n}\left\Vert
p_{0}-p_{\min ,\lambda \left( \delta \right) ,\alpha \left( \delta \right)
}\right\Vert _{H^{2}\left( 0,1\right) }, \, n=1,2,...,  \label{4.8}
\end{equation}%
\emph{\ }%
\begin{equation}
\left\Vert y^{n}-y_{\min ,\lambda \left( \delta \right) ,\alpha \left(
\delta \right) }\right\Vert _{H^{2}\left( 0,1\right) }\leq q^{n}\left\Vert
y^{0}-y_{\min ,\lambda \left( \delta \right) ,\alpha \left( \delta \right)
}\right\Vert _{H^{2}\left( 0,1\right) }, \, n=1,2,...,  \label{4.9}
\end{equation}%
\begin{equation}
\left\Vert y^{n}-y^{\ast }\right\Vert _{H^{2}\left( 0,1\right) }\leq
C_{2}\delta ^{1/4}+q^{n}\left\Vert y^{0}-y_{\min ,\lambda \left( \delta
\right) ,\alpha \left( \delta \right) }\right\Vert _{H^{2}\left( 0,1\right)
}, \, n=1,2,...,  \label{4.10}
\end{equation}%
\begin{equation}
\left\Vert c^{n}-c^{\ast }\right\Vert _{L_{2}\left( 0,1\right) }\leq
C_{2}\delta ^{1/4}+q^{n}\left\Vert y^{0}-y_{\min ,\lambda \left( \delta
\right) ,\alpha \left( \delta \right) }\right\Vert _{H^{2}\left( 0,1\right)
}, \, n=1,2,...,  \label{4.11}
\end{equation}%
\emph{where functions }$c^{n}\left( x\right) $\emph{\ are found from
components }$y_{m}^{n},m=0,...,N-1$\emph{\ of vector functions }$y^{n}$\emph{%
\ via formula (\ref{3.8}).}

\textbf{Remark 4.1}. \emph{Since }$R>0$\emph{\ is an arbitrary number and }$%
p_{0}$ \emph{is an arbitrary point of the ball }$B_{0}\left( R\right) $\emph{%
, then Theorem 4.5 claims the global convergence of the gradient projection
method (\ref{4.7}), see section 1 for our definition of this term. }

\section{Proofs}

\label{sec:5}

In this section we prove theorems formulated in section 4. The proof of
Corollary 4.1 is completely similarly with the proof of Theorem 4.1. The
same is true about the similarity of the proof of Theorem 4.2 with the proof
of Theorem 3.1 of \cite{BakKlib}. Hence, we omit proofs of Corollary 4.1 and
Theorem 4.2.

\subsection{Proof of Theorem 4.1}

\label{sec:5.1}

Below $\left( ,\right) _{2}$ is the scalar product in $H^{2}\left(
0,1\right) .$ Since we work with complex valued $N-$D vector functions, it
is convenient to consider them as pairs of real valued functions, 
\[
p\left( x\right) =\left( \func{Re}p\left( x\right) ,\func{Im}p\left(
x\right) \right) =\left( p_{1}\left( x\right) ,p_{2}\left( x\right) \right) .
\]%
Thus, in fact we work with real valued $N-$D vector functions $p_{1}\left(
x\right) ,p_{2}\left( x\right) $ and the vector function $p\left( x\right) $
is $2N-$D. We use standard definitions of scalar products for Hilbert spaces
of vector functions. Also, for brevity, notations of those spaces of vector
functions are the same as ones for regular functions. The use of the complex
conjugation below is for convenience only. Below we use the following formula%
\begin{equation}
\left\vert a\right\vert ^{2}-\left\vert b\right\vert ^{2}=\left( a-b\right) 
\overline{a}+\left( \overline{a}-\overline{b}\right) b, \quad \forall
a,b\in \mathbb{C}.  \label{5.0}
\end{equation}%
For complex valued vectors $a=\left( a_{1},a_{2}\right) ,b=\left(
b_{1},b_{2}\right) $ with $a_{1},a_{2},b_{1},b_{2}\in \mathbb{R}^{N}$ and
with the notation $\left[ ,\right] $ for the scalar product in $\mathbb{R}%
^{2N}$ (\ref{5.0}) becomes%
\begin{equation}
\left\vert a\right\vert ^{2}-\left\vert b\right\vert ^{2}=\left[ a-b,%
\overline{a}\right] +\left[ \overline{a}-\overline{b},b\right] .
\label{5.01}
\end{equation}

Let vector functions $y_{\left( 1\right) },y_{\left( 2\right) }\in B\left(
2R,f_{0},f_{1}\right) .$ Denote $h=y_{\left( 2\right) }-y_{\left( 1\right) }.
$ By (\ref{3.13}) 
\begin{equation}
h\left( 0\right) =h^{\prime }\left( 0\right) =h^{\prime }\left( 1\right) =0,
\label{5.1}
\end{equation}%
\begin{equation}
\left\Vert h\right\Vert _{H^{2}\left( 0,1\right) }\leq 4R.  \label{5.2}
\end{equation}%
In particular, (\ref{5.1}) implies that $h\in H_{0}^{2}\left( 0,1\right) .$
Also, by embedding theorem $H^{2}\left( 0,1\right) \subset C^{1}\left[ 0,1%
\right] .$ Hence, (\ref{5.2}) implies that with a generic constant $C$, 
\begin{equation}
\left\Vert h^{\prime }\right\Vert _{C\left[ 0,1\right] }\leq CR.
\label{5.20}
\end{equation}%
Keeping in mind of using (\ref{5.01}), denote 
\begin{equation}
a=\left( y_{\left( 1\right) }+h\right) ^{\prime \prime }+F\left( y_{\left(
1\right) }^{\prime }+h^{\prime }\right) ,\, b=y_{\left( 1\right)
}^{\prime \prime }+F\left( y_{\left( 1\right) }^{\prime }\right) .
\label{5.3}
\end{equation}%
Hence, 
\[
a-b=h^{\prime \prime }+F\left( y_{\left( 1\right) }^{\prime }+h^{\prime
}\right) -F\left( y_{\left( 1\right) }^{\prime }\right) .
\]%
Recalling that the vector function $F\left( y^{\prime }\right) $ is
quadratic with respect to the components $y_{n}^{\prime }\left( x\right) $
of the vector function $y^{\prime },$ we obtain%
\begin{equation}
F\left( y_{\left( 1\right) }^{\prime }+h^{\prime }\right) -F\left( y_{\left(
1\right) }^{\prime }\right) =G_{1}\left( y_{\left( 1\right) }^{\prime
},h^{\prime }\right) +G_{2}\left( y_{\left( 1\right) }^{\prime },h^{\prime
}\right) ,  \label{5.4}
\end{equation}%
\begin{equation}
a-b=h^{\prime \prime }+G_{1}\left( y_{\left( 1\right) }^{\prime },h^{\prime
}\right) +G_{2}\left( y_{\left( 1\right) }^{\prime },h^{\prime }\right) ,
\label{5.5}
\end{equation}%
where the vector function $G_{1}\left( y_{\left( 1\right) }^{\prime
},h^{\prime }\right) $ is linear with respect to each $y_{\left( 1\right)
}^{\prime }$ and $h^{\prime }.$ The following estimates are valid: 
\begin{equation}
\left\vert G_{1}\left( y_{\left( 1\right) }^{\prime },h^{\prime }\right)
\right\vert \leq C_{1}\left\vert h^{\prime }\right\vert ,  \label{5.50}
\end{equation}%
\begin{equation}
\left\vert G_{2}\left( y_{\left( 1\right) }^{\prime },h^{\prime }\right)
\right\vert \leq C_{1}\left\vert h^{\prime }\right\vert ^{2}.  \label{5.6}
\end{equation}%
Hence, keeping in mind (\ref{5.01}), (\ref{5.3})-(\ref{5.5}), we obtain 
\[
\left[ a-b,\overline{a}\right] =
\]%
\[
\left[ h^{\prime \prime }+G_{1}\left( y_{\left( 1\right) }^{\prime
},h^{\prime }\right) +G_{2}\left( y_{\left( 1\right) }^{\prime },h^{\prime
}\right) ,\overline{\left( y_{\left( 1\right) }+h\right) ^{\prime \prime
}+F\left( y_{\left( 1\right) }^{\prime }+h^{\prime }\right) }\right] 
\]%
\begin{equation}
=\left[ h^{\prime \prime }+G_{1}\left( y_{\left( 1\right) }^{\prime
},h^{\prime }\right) ,\overline{y_{\left( 1\right) }^{\prime \prime
}+F\left( y_{\left( 1\right) }^{\prime }\right) }\right] +\left\vert h^{\prime \prime }\right\vert ^{2}
 \label{5.7}
\end{equation}%
\[
+\left[ h^{\prime \prime },%
\overline{G_{1}\left( y_{\left( 1\right) }^{\prime },h^{\prime }\right)
+G_{2}\left( y_{\left( 1\right) }^{\prime },h^{\prime }\right) }\right] +%
\left[ G_{1}\left( y_{\left( 1\right) }^{\prime },h^{\prime }\right)
+G_{2}\left( y_{\left( 1\right) }^{\prime },h^{\prime }\right) ,\overline{%
h^{\prime \prime }}\right]  
\]
\[
+\left[ G_{1}\left( y_{\left( 1\right) }^{\prime },h^{\prime }\right) ,%
\overline{G_{2}\left( y_{\left( 1\right) }^{\prime },h^{\prime }\right) }%
\right] +\left[ G_{2}\left( y_{\left( 1\right) }^{\prime },h^{\prime
}\right) ,\overline{G_{1}\left( y_{\left( 1\right) }^{\prime },h^{\prime
}\right) }\right] 
\]%
\[
+\left\vert G_{1}\left( y_{\left( 1\right) }^{\prime },h^{\prime }\right)
\right\vert ^{2}+\left\vert G_{2}\left( y_{\left( 1\right) }^{\prime
},h^{\prime }\right) \right\vert ^{2}+\left[ G_{2}\left( y_{\left( 1\right)
}^{\prime },h^{\prime }\right) ,\overline{y_{\left( 1\right) }^{\prime
\prime }+F\left( y_{\left( 1\right) }^{\prime }\right) }\right] 
\]%
Let $S_{1}\left( x\right) $ denotes the sum of 3rd, 4th and 5th lines of (%
\ref{5.7}),%
\[
S_{1}\left( x\right) =\left\vert h^{\prime \prime }\right\vert ^{2}+\left[
h^{\prime \prime },\overline{G_{1}\left( y_{\left( 1\right) }^{\prime
},h^{\prime }\right) +G_{2}\left( y_{\left( 1\right) }^{\prime },h^{\prime
}\right) }\right] 
\]
\[
+\left[ G_{1}\left( y_{\left( 1\right) }^{\prime
},h^{\prime }\right) +G_{2}\left( y_{\left( 1\right) }^{\prime },h^{\prime
}\right) ,\overline{h^{\prime \prime }}\right] +\left[ G_{1}\left( y_{\left( 1\right) }^{\prime },h^{\prime }\right) ,%
\overline{G_{2}\left( y_{\left( 1\right) }^{\prime },h^{\prime }\right) }%
\right]
\]
\begin{equation}
 +\left[ G_{2}\left( y_{\left( 1\right) }^{\prime },h^{\prime
}\right) ,\overline{G_{1}\left( y_{\left( 1\right) }^{\prime },h^{\prime
}\right) }\right] +\left\vert G_{1}\left( y_{\left( 1\right) }^{\prime },h^{\prime }\right)
\right\vert ^{2}+\left\vert G_{2}\left( y_{\left( 1\right) }^{\prime
},h^{\prime }\right) \right\vert ^{2}
  \label{5.70}
\end{equation}%
\[
+\left[ G_{2}\left( y_{\left( 1\right)
}^{\prime },h^{\prime }\right) ,\overline{y_{\left( 1\right) }^{\prime
\prime }+F\left( y_{\left( 1\right) }^{\prime }\right) }\right] .
\]%
Then, using the inequality 
\[
\left[ d_{1},d_{2}\right] \geq -\frac{1}{4}\left\vert d_{1}\right\vert
^{2}-\left\vert d_{2}\right\vert ^{2}, \quad \forall d_{1},d_{2}\in \mathbb{R}^{2N},
\]%
(\ref{5.20}), (\ref{5.50}) and (\ref{5.6}), we obtain the following estimate
from the below for $S_{1}\left( x\right) $, 
\begin{equation}
S_{1}\left( x\right) \geq \frac{1}{2}\left\vert h^{\prime \prime }\left(
x\right) \right\vert ^{2}-C_{1}\left\vert h^{\prime }\left( x\right)
\right\vert ^{2}, \quad  x\in \left( 0,1\right) .  \label{5.8}
\end{equation}%
By (\ref{5.7}) $\left[ a-b,\overline{a}\right] $ can be written as%
\begin{equation}
\left[ a-b,\overline{a}\right] =\left[ h^{\prime \prime }+G_{1}\left(
y_{\left( 1\right) }^{\prime },h^{\prime }\right) ,\overline{y_{\left(
1\right) }^{\prime \prime }+F\left( y_{\left( 1\right) }^{\prime }\right) }%
\right] +S_{1}.  \label{5.9}
\end{equation}

Similarly the second term in the right hand side of (\ref{5.01}) can be
represented as%
\[
\left[ \overline{a}-\overline{b},b\right] =\left[ \overline{h^{\prime \prime
}+G_{1}\left( y_{\left( 1\right) }^{\prime },h^{\prime }\right) +G_{2}\left(
y_{\left( 1\right) }^{\prime },h^{\prime }\right) },y_{\left( 1\right)
}^{\prime \prime }+F\left( y_{\left( 1\right) }^{\prime }\right) \right] 
\]%
\begin{equation}
=\left[ \overline{h^{\prime \prime }+G_{1}\left( y_{\left( 1\right)
}^{\prime },h^{\prime }\right) },y_{\left( 1\right) }^{\prime \prime
}+F\left( y_{\left( 1\right) }^{\prime }\right) \right]   \label{5.10}
\end{equation}%
\[
+\left[ \overline{G_{2}\left( y_{\left( 1\right) }^{\prime },h^{\prime
}\right) },y_{\left( 1\right) }^{\prime \prime }+F\left( y_{\left( 1\right)
}^{\prime }\right) \right] ,
\]%
where by (\ref{5.6})%
\begin{equation}
\left[ \overline{G_{2}\left( y_{\left( 1\right) }^{\prime },h^{\prime
}\right) },y_{\left( 1\right) }^{\prime \prime }+F\left( y_{\left( 1\right)
}^{\prime }\right) \right] =S_{2}\left( x\right) \geq -C_{1}\left\vert
h^{\prime }\right\vert ^{2}.  \label{5.11}
\end{equation}%
Thus, (\ref{5.0}), (\ref{5.3}) and (\ref{5.7})-(\ref{5.11}) imply that 
\[
\left\vert a\right\vert ^{2}-\left\vert b\right\vert ^{2}=Lin\left( h\right)
+S\left( x\right) ,
\]%
where the term $Lin\left( h\right) $ is linear with respect to the vector
function $h=\left( h_{1},h_{2}\right) $, 
\begin{equation}
Lin\left( h\right) =\left[ h^{\prime \prime }+G_{1}\left( y_{\left( 1\right)
}^{\prime },h^{\prime }\right) ,\overline{y_{\left( 1\right) }^{\prime
\prime }+F\left( y_{\left( 1\right) }^{\prime }\right) }\right] 
\label{5.12}
\end{equation}%
\[
+\left[ \overline{h^{\prime \prime }+G_{1}\left( y_{\left( 1\right)
}^{\prime },h^{\prime }\right) },y_{1}^{\prime \prime }+F\left( y_{\left(
1\right) }^{\prime }\right) \right] 
\]%
and the function $S\left( x\right) $ is%
\begin{equation}
S\left( x\right) =S_{1}\left( x\right) +S_{2}\left( x\right) \geq \frac{1}{2}%
\left\vert h^{\prime \prime }\left( x\right) \right\vert
^{2}-C_{1}\left\vert h^{\prime }\left( x\right) \right\vert ^{2}.
\label{5.13}
\end{equation}%
Next, it follows from (\ref{5.6}), (\ref{5.7}), (\ref{5.70}) and (\ref{5.11}%
) that 
\begin{equation}
\left\vert S\left( x\right) \right\vert \leq C_{1}\left( \left\vert
h^{\prime \prime }\right\vert ^{2}+\left\vert h^{\prime }\right\vert
^{2}\right) \left( x\right) .  \label{5.16}
\end{equation}

Hence, (\ref{3.14}) implies that%
\[
J_{\lambda ,\alpha }\left( y_{\left( 1\right) }+h\right) -J_{\lambda ,\alpha
}\left( y_{\left( 1\right) }\right) 
\]%
\begin{equation}
=e^{2\lambda }\int\limits_{0}^{1}Lin\left( h\right) e^{-2\lambda x}dx+\alpha
\left( h,y_{\left( 1\right) }\right) _{2}+\alpha \left( y_{\left( 1\right)
},h\right) _{2}  \label{5.14}
\end{equation}%
\[
+e^{2\lambda }\int\limits_{0}^{1}S\left( x\right) e^{-2\lambda x}dx+\alpha
\left\Vert h\right\Vert _{H^{2}\left( 0,1\right) }^{2}.
\]%
Let $L_{\lambda ,\alpha }\left( h\right) $ denotes the second line of (\ref%
{5.14}). Then (\ref{5.1}), (\ref{5.50}) and (\ref{5.12}) imply that $%
L_{\lambda ,\alpha }\left( h\right) :H_{0}^{2}\left( 0,1\right) \rightarrow 
\mathbb{R}$ is a bounded linear functional. Hence, by Riesz theorem, there
exists an element $w_{\lambda ,\alpha }\in H_{0}^{2}\left( 0,1\right) $ such
that 
\begin{equation}
L_{\lambda ,\alpha }\left( h\right) =\left( w_{\lambda ,\alpha },h\right)
_{2}, \quad  \forall h=\left( h_{1},h_{2}\right) \in H_{0}^{2}\left( 0,1\right) .
\label{5.15}
\end{equation}%
Hence, (\ref{5.16}), (\ref{5.14}) and (\ref{5.15}) imply that%
\begin{equation}
J_{\lambda ,\alpha }\left( y_{\left( 1\right) }+h\right) -J_{\lambda ,\alpha
}\left( y_{\left( 1\right) }\right) =\left( w_{\lambda ,\alpha },h\right)
_{2}+o\left( \left\Vert h\right\Vert _{H^{2}\left( 0,1\right) }\right) ,
\label{5.17}
\end{equation}%
as $\left\Vert h\right\Vert _{H^{2}\left( 0,1\right) }\rightarrow 0.$ Hence, 
$w_{\lambda ,\alpha }\in H_{0}^{2}\left( 0,1\right) $ is the Frech\'{e}t
derivative $J_{\lambda ,\alpha }^{\prime }\left( y_{\left( 1\right) }\right) 
$ of the functional $J_{\lambda ,\alpha }$ at the point $y_{\left( 1\right)
},$ i.e. $J_{\lambda ,\alpha }^{\prime }\left( y_{\left( 1\right) }\right)
=w_{\lambda ,\alpha }.$ Thus, the existence of the Frech\'{e}t derivative of
the functional $J_{\lambda ,\alpha }\left( y\right) $ at any point $y\in
B\left( 2R,f_{0},f_{1}\right) $ is established.

Below we prove the strict convexity estimate (\ref{4.1}). Let now $y_{\left(
1\right) }\in \overline{B\left( R,f_{0},f_{1}\right) }$ and $y_{\left(
2\right) }\in \overline{B\left( R,f_{0},f_{1}\right) }$ be two arbitrary
points. We keep the same notation for their difference $h=y_{\left( 1\right)
}-y_{\left( 2\right) }.$ The vector function $h$ satisfies boundary
conditions (\ref{5.1}), and in (\ref{5.2}) $4R$ should be replaced with $2R.$
Thus, by (\ref{5.13}) and (\ref{5.14}) 
\[
J_{\lambda ,\alpha }\left( y_{\left( 1\right) }+h\right) -J_{\lambda ,\alpha
}\left( y_{\left( 1\right) }\right) -J_{\lambda ,\alpha }^{\prime }\left(
y_{\left( 1\right) }\right) \left( h\right) 
\]%
\begin{equation}
=e^{2\lambda}\int\limits_{0}^{1}S\left( x\right) e^{-2\lambda x}dx +\alpha \left\Vert
h\right\Vert _{H^{2}\left( 0,1\right) }^{2} \geq \frac{1}{2}e^{2\lambda }\int\limits_{0}^{1}\left\vert h^{\prime \prime
}\right\vert ^{2}e^{-2\lambda x}dx
\label{5.18}
\end{equation}
\[
-C_{1}e^{2\lambda
}\int\limits_{0}^{1}\left\vert h^{\prime }\right\vert ^{2}e^{-2\lambda
x}dx+\alpha \left\Vert h\right\Vert _{H^{2}\left( 0,1\right) }^{2}.
\]

Using Lemma 3.2 and boundary conditions (\ref{5.1}), we estimate from the
below the second line of (\ref{5.18}) as 
\[
\frac{1}{2}e^{2\lambda }\int\limits_{0}^{1}\left\vert h^{\prime \prime
}\right\vert ^{2}e^{-2\lambda x}dx-C_{1}e^{2\lambda
}\int\limits_{0}^{1}\left\vert h^{\prime }\right\vert ^{2}e^{-2\lambda
x}dx+\alpha \left\Vert h\right\Vert _{H^{2}\left( 0,1\right) }^{2}
\]%
\begin{equation}
\geq Ce^{2\lambda }\int\limits_{0}^{1}\left\vert h^{\prime \prime
}\right\vert ^{2}e^{-2\lambda x}dx+C\lambda e^{2\lambda
}\int\limits_{0}^{1}|h^{\prime }|^{2}e^{-2\lambda x}dx \label{5.180}
\end{equation}%
\[
+C\lambda^{3}e^{2\lambda }\int\limits_{0}^{1}|h|^{2}e^{-2\lambda x}dx  -C_{1}e^{2\lambda }\int\limits_{0}^{1}\left\vert h^{\prime }\right\vert
^{2}e^{-2\lambda x}dx+\alpha \left\Vert h\right\Vert _{H^{2}\left(
0,1\right) }^{2},\forall \lambda \geq 1.
\]%
Hence, one can choose a sufficiently large number $\lambda _{1}=\lambda
_{1}\left( R,F,N\right) >1$ such that for all $\lambda \geq \lambda _{1}$ we
have $C_{1}<C\lambda /2.$ Note that $e^{2\lambda }e^{-2\lambda x}>1$ for $%
x\in \left( 0,1\right) .$ Hence, (\ref{5.18}) and (\ref{5.180}) imply%
\[
J_{\lambda ,\alpha }\left( y_{\left( 1\right) }+h\right) -J_{\lambda ,\alpha
}\left( y_{\left( 1\right) }\right) -J_{\lambda ,\alpha }^{\prime }\left(
y_{\left( 1\right) }\right) \left( h\right) 
\]
\[
\geq \frac{1}{2}e^{2\lambda }\int\limits_{0}^{1}\left\vert h^{\prime \prime
}\right\vert ^{2}e^{-2\lambda x}dx-C_{1}e^{2\lambda
}\int\limits_{0}^{1}\left\vert h^{\prime }\right\vert ^{2}e^{-2\lambda
x}dx+\alpha \left\Vert h\right\Vert _{H^{2}\left( 0,1\right) }^{2}
\]%
\[
\geq Ce^{2\lambda }\int\limits_{0}^{1}\left\vert h^{\prime \prime
}\right\vert ^{2}e^{-2\lambda x}dx+C_{1}\lambda e^{2\lambda
}\int\limits_{0}^{1}|h^{\prime }|^{2}e^{-2\lambda x}dx
\]%
\[
+C_{1}\lambda ^{3}e^{2\lambda }\int\limits_{0}^{1}|h|^{2}e^{-2\lambda
x}dx+\alpha \left\Vert h\right\Vert _{H^{2}\left( 0,1\right) }^{2}
\]%
\[
\geq C_{1}\int\limits_{0}^{1}\left( \left\vert h^{\prime \prime }\right\vert
^{2}+|h^{\prime }|^{2}+|h|^{2}\right) dx+\alpha \left\Vert h\right\Vert
_{H^{2}\left( 0,1\right) }^{2}
\]%
\[
=\left( C_{1}+\alpha \right) \left\Vert h\right\Vert _{H^{2}\left(
0,1\right) }^{2}\geq C_{1}\left\Vert h\right\Vert _{H^{2}\left( 0,1\right)
}^{2},
\]%
which is equivalent with (\ref{4.1}). $\square $

It is important for the proof of Theorem 4.4 to point out that, until (\ref%
{5.15}), we have not used boundary conditions (\ref{5.1}) for the vector
function $h\left( x\right) $. Hence, (\ref{5.50}), (\ref{5.12}), (\ref{5.16}%
) and (\ref{5.14}) imply that 
\begin{equation}
\left\vert J_{\lambda ,\alpha }\left( y+z\right) -J_{\lambda ,\alpha }\left(
y\right) \right\vert \leq C_{1}\left\Vert z\right\Vert _{H^{2}\left(
0,1\right) }e^{2\lambda },  \label{5.181}
\end{equation}%
\begin{equation}
\forall y\in B\left( 2R,f_{0},f_{1}\right) ,\quad \forall z\in \left\{ z\in
H^{2}\left( 0,1\right) :\left\Vert z\right\Vert _{H^{2}\left( 0,1\right)
}<4R\right\} .  \label{5.182}
\end{equation}

\subsection{Proof of Theorem 4.3}

\label{sec:5.2}

Let $f\left( x\right) \in C^{2}\left[ 0,1\right] $ be the function defined
in (\ref{1}), (\ref{2}). Then Proposition 3.1 implies that 
\begin{equation}
p=y-f\in \overline{B_{0}\left( R\right) }, \quad \forall y\in \overline{B\left(
R,f_{0},f_{1}\right) }.  \label{5.21}
\end{equation}%
It follows immediately from Theorem 2.1 of \cite{BakKlib} and Corollary 4.1
that for any $\lambda \geq \lambda _{2}$ there exists unique minimizer $%
p_{\min ,\lambda ,\alpha }$ of the functional $\Phi _{\lambda ,\alpha
}\left( p\right) $ on the set $\overline{B_{0}\left( R\right) },$%
\begin{equation}
\min_{p\in \overline{B_{0}\left( R\right) }}\Phi _{\lambda ,\alpha }\left(
p\right) =\Phi _{\lambda ,\alpha }\left( p_{\min ,\lambda ,\alpha }\right) .
\label{5.210}
\end{equation}

Consider now the vector function $\widetilde{y}_{\lambda ,\alpha },$%
\begin{equation}
\widetilde{y}_{\lambda ,\alpha }\left( x\right) =p_{\min ,\lambda ,\alpha
}\left( x\right) +f\left( x\right) .  \label{5.23}
\end{equation}%
By Proposition 3.1 the function $\widetilde{y}_{\lambda ,\alpha }$ $\in 
\overline{B\left( R,f_{0},f_{1}\right) }.$ Let $y\neq \widetilde{y}_{\lambda
,\alpha }$ be an arbitrary vector function from $\overline{B\left(
R,f_{0},f_{1}\right) }.$ Denote $p_{y}=y-f.$ Then by Proposition 3.1 $%
p_{y}\in \overline{B_{0}\left( R\right) }.$ We have: 
\begin{equation}
J_{\lambda ,\alpha }\left( y\right) =J_{\lambda ,\alpha }\left( \left(
y-f\right) +f\right) =J_{\lambda ,\alpha }\left( p_{y}+f\right) .
\label{5.24}
\end{equation}%
And by (\ref{5.22}) 
\begin{equation}
J_{\lambda ,\alpha }\left( p_{y}+f\right) =\Phi _{\lambda ,\alpha }\left(
p_{y}\right) .  \label{5.240}
\end{equation}%
Since $p_{y}\neq p_{\min ,\lambda ,\alpha },$ then 
\begin{equation}
\Phi _{\lambda ,\alpha }\left( p_{y}\right) >\Phi _{\lambda ,\alpha }\left(
p_{\min ,\lambda ,\alpha }\right) .  \label{5.25}
\end{equation}%
On the other hand, by (\ref{5.22}) and (\ref{5.23})-(\ref{5.25}) 
\begin{equation}
J_{\lambda ,\alpha }\left( y\right) =\Phi _{\lambda ,\alpha }\left(
p_{y}\right) >\Phi _{\lambda ,\alpha }\left( p_{\min ,\lambda ,\alpha
}\right) =J_{\lambda ,\alpha }\left( \widetilde{y}_{\lambda ,\alpha }\right)
.  \label{5.26}
\end{equation}%
Hence, (\ref{5.26}) implies that $\widetilde{y}_{\lambda ,\alpha }$ is the
unique minimizer $y_{\min ,\lambda ,\alpha }$ of the functional $J_{\lambda
,\alpha }$ on the set $\overline{B\left( R,f_{0},f_{1}\right) },$%
\begin{equation}
\min_{\overline{B\left( R,f_{0},f_{1}\right) }}J_{\lambda ,\alpha }\left(
y\right) =J_{\lambda ,\alpha }\left( y_{\min ,\lambda ,\alpha }\right) .
\label{5.27}
\end{equation}%
Inequalities (\ref{4.20}) and (\ref{4.2}) follow immediately from Lemma 2.1
of \cite{BakKlib}, (\ref{5.210}) and (\ref{5.27}) respectively. $\square $

\subsection{Proof of Theorem 4.4}

\label{sec:5.3}

Since $\left( y^{\ast }\right) ^{\prime \prime }+F\left( \left( y^{\ast
}\right) ^{\prime }\right) =0$ and $\alpha =\sqrt{\delta },$ then (\ref{3.14}%
) implies that%
\begin{equation}
J_{\lambda ,\alpha }\left( y^{\ast }\right) =\alpha \left\Vert y^{\ast
}\right\Vert _{H^{2}\left( 0,1\right) }^{2}\leq \sqrt{\delta }R^{2}.
\label{5.28}
\end{equation}%
Recall that 
\begin{equation}
y^{\ast }-f^{\ast }=p^{\ast }\in B_{0}\left( R\right) .  \label{5.280}
\end{equation}%
We have%
\begin{equation}
J_{\lambda ,\alpha }\left( p^{\ast }+f\right) =J_{\lambda ,\alpha }\left(
\left( p^{\ast }+f^{\ast }\right) +\left( f-f^{\ast }\right) \right)
=J_{\lambda ,\alpha }\left( y^{\ast }+\left( f-f^{\ast }\right) \right) .
\label{5.29}
\end{equation}%
Temporary denote $h=f-f^{\ast }.$ Since by Proposition 3.1 $\left( p^{\ast
}+f\right) \in B\left( R,f_{0},f_{1}\right) $, then, using (\ref{3.131}), (%
\ref{5.181}), (\ref{5.182}), (\ref{5.28}) and (\ref{5.280}), we obtain%
\begin{equation}
J_{\lambda ,\alpha }\left( p^{\ast }+f\right) =J_{\lambda ,\alpha }\left(
y^{\ast }+h\right) \leq J_{\lambda ,\alpha }\left( y^{\ast }\right) +B\delta
e^{2\lambda }\leq C_{2}\left( \sqrt{\delta }+\delta e^{2\lambda }\right) ,
\label{5.30}
\end{equation}%
see the line below (\ref{3.131}) for the definition of the number $B=B\left(
\chi \right) >0.$ Choose a number $\delta _{0}\in \left( 0,e^{-4\lambda
_{2}}\right) .$ Let $\delta \in \left( 0,\delta _{0}\right) .$ Choose $%
\lambda =\lambda \left( \delta \right) =\ln \left( \delta ^{-1/4}\right) .$
Then $\lambda \left( \delta \right) >\lambda _{2}$ and $\delta e^{2\lambda }=%
\sqrt{\delta }$. Hence, (\ref{5.30}) implies that%
\begin{equation}
J_{\lambda ,\alpha }\left( p^{\ast }+f\right) \leq C_{2}\sqrt{\delta }.
\label{5.300}
\end{equation}

Let 
\begin{equation}
p_{\min ,\lambda ,\alpha }=y_{\min ,\lambda ,\alpha }-f.  \label{5.31}
\end{equation}%
By (\ref{5.22})%
\[
J_{\lambda ,\alpha }\left( p^{\ast }+f\right) -J_{\lambda ,\alpha }\left(
p_{\min ,\lambda ,\alpha }+f\right) =\Phi _{\lambda ,\alpha }\left( p^{\ast
}\right) -\Phi _{\lambda ,\alpha }\left( p_{\min ,\lambda ,\alpha }\right) . 
\]%
Hence, by Corollary 4.1,%
\[
J_{\lambda ,\alpha }\left( p^{\ast }+f\right) -J_{\lambda ,\alpha }\left(
p_{\min ,\lambda ,\alpha }+f\right) 
\]
\[
-J_{\lambda ,\alpha }^{\prime }\left(
p_{\min ,\lambda ,\alpha }+f\right) \left( \left( p^{\ast }+f\right) -\left(
p_{\min ,\lambda ,\alpha }+f\right) \right) 
\]
\begin{equation}
=\Phi _{\lambda ,\alpha }\left( p^{\ast }\right) -\Phi _{\lambda ,\alpha
}\left( p_{\min ,\lambda ,\alpha }\right) -\Phi _{\lambda ,\alpha }^{\prime
}\left( p_{\min ,\lambda ,\alpha }\right) \left( p^{\ast }-p_{\min ,\lambda
,\alpha }\right)  \label{5.32}
\end{equation}%
\[
\geq C_{2}\left\Vert p^{\ast }-p_{\min ,\lambda ,\alpha }\right\Vert
_{H^{2}\left( 0,1\right) }^{2}. 
\]%
Next, by (\ref{4.20}) $-\Phi _{\lambda ,\alpha }^{\prime }\left( p_{\min
,\lambda ,\alpha }\right) \left( p^{\ast }-p_{\min ,\lambda ,\alpha }\right)
\leq 0.$ Hence, (\ref{5.32}) implies that 
\begin{equation}
\Phi _{\lambda ,\alpha }\left( p^{\ast }\right) \geq C_{2}\left\Vert p^{\ast
}-p_{\min ,\lambda ,\alpha }\right\Vert _{H^{2}\left( 0,1\right) }^{2}.
\label{5.320}
\end{equation}%
On the other hand, since by (\ref{5.22}) $\Phi _{\lambda ,\alpha }\left(
p^{\ast }\right) =J_{\lambda ,\alpha }\left( p^{\ast }+f\right) ,$ then,
using (\ref{5.300}) and (\ref{5.320}), we obtain 
\[
\left\Vert p^{\ast }-p_{\min ,\lambda ,\alpha }\right\Vert _{H^{2}\left(
0,1\right) }\leq C_{2}\delta ^{1/4}, 
\]%
which establishes (\ref{4.5}). Next, 
\[
\left\Vert p^{\ast }-p_{\min ,\lambda ,\alpha }\right\Vert _{H^{2}\left(
0,1\right) }=\left\Vert \left( p^{\ast }+f^{\ast }\right) -\left( p_{\min
,\lambda ,\alpha }+f\right) +\left( f-f^{\ast }\right) \right\Vert
_{H^{2}\left( 0,1\right) } 
\]%
\begin{equation}
\geq \left\Vert y^{\ast }-y_{\min ,\lambda ,\alpha }\right\Vert
_{H^{2}\left( 0,1\right) }-\left\Vert f-f^{\ast }\right\Vert _{H^{2}\left(
0,1\right) }  \label{5.34}
\end{equation}%
\[
\geq \left\Vert y^{\ast }-y_{\min ,\lambda ,\alpha }\right\Vert
_{H^{2}\left( 0,1\right) }-B\delta . 
\]%
Here, to obtain the term $B\delta $, we have used (\ref{3.131}). Hence,
using (\ref{4.5}) and (\ref{5.34}), we obtain $\left\Vert y^{\ast }-y_{\min
,\lambda ,\alpha }\right\Vert _{H^{2}\left( 0,1\right) }\leq C_{2}\delta
^{1/4}$, which is the same as estimate (\ref{4.50}). Estimate (\ref{4.6})
follows immediately from (\ref{3.8}) and (\ref{4.50}). \ $\square $

\subsection{Proof of Theorem 4.5}

\label{sec:5.4}

Estimate (\ref{4.8}) follows immediately from the combination of Corollary
4.1 and Theorem 4.3 with Theorem 2.1 of \cite{BakKlib}. Given that $%
y^{n}=p_{n}+f$ and also that by Theorem 4.3 $y_{\min ,\lambda ,\alpha
}=p_{\min ,\lambda ,\alpha }+f$ , (\ref{4.8}) implies (\ref{4.9}). Next, by
triangle inequality%
\begin{equation}
\left\Vert y^{n}-y^{\ast }\right\Vert _{H^{2}\left( 0,1\right) }\leq
\left\Vert y_{\min ,\lambda \left( \delta \right) ,\alpha \left( \delta
\right) }-y^{\ast }\right\Vert _{H^{2}\left( 0,1\right) }+\left\Vert
y^{n}-y_{\min ,\lambda \left( \delta \right) ,\alpha \left( \delta \right)
}\right\Vert _{H^{2}\left( 0,1\right) }.  \label{5.35}
\end{equation}%
Hence, estimate (\ref{4.10}) follows from (\ref{4.50}), (\ref{4.9}) and (\ref%
{5.35}). Finally, (\ref{4.11}) follows from (\ref{3.8}), (\ref{4.50}) and (%
\ref{4.10}). $\square $

\section{Numerical studies}

\label{sec:6}

In this section we present a numerical verification of proposed method by
reconstructing the coefficient $c(x)$ from both computationally simulated
and experimental data. Since in our experiments we have to image only a
single target for each data set, then we focus in our computationally
simulated data on the case of a single inclusion only.

\subsection{Data generation}

\label{sec:6.1}

To generate computationally simulated data for our CIP, we solve the forward
problem (\ref{2.4}), (\ref{2.6}) via solving the one dimensional
Lippmann-Schwinger equation (\ref{2.192}) in which the function $c\left(
x\right) $ is set as follows: 
\begin{equation}
c\left( x\right) :=c_{true}(x)=\left\{ 
\begin{array}{cc}
\widehat{c}_{true}, & \quad x\in (x_{loc}-d/2,x_{loc}+d/2), \\ 
1, & \quad \mbox{elsewhere}.%
\end{array}%
\right.  \label{eq:c_true}
\end{equation}%
Here, $\widehat{c}_{true}=const.>0$ is the dielectric constant of a
simulated target, $x_{loc}$ is the location of the center of this target,
and $d$ is its width. Thus, targets used in computational simulations are
step-wise functions. In numerical experiments we consider the following sets
of dielectric constants $\widehat{c}_{true}$ and locations $x_{loc}$: 
\begin{equation}
\widehat{c}_{true}=\{3.0,\,4.0,\,5.0,\,6.0\},\quad
x_{loc}=\{0.1,\,0.2,\,0.3,\,0.4\}.  \label{eq:targers}
\end{equation}%
For this study we use the same width $d=0.1$ in (\ref{eq:c_true}) for all
targets. Since we test four locations of centers for each of these four
values of $\widehat{c}_{true},$ then we reconstruct total sixteen (16)
targets.

We generate our simulated data for the interval of wave numbers $k_{m}\in
\lbrack 0.5,\,1.5]$. This interval is divided into $N_{k}$ equal
subintervals. Hence, we obtain a grid of equally spaced points: $k_{m}=%
\underline{k}+mh_{k}$, $m=0,\dots ,N_{k}$, where $h_{k}=(\overline{k}-%
\underline{k})/N_{k}.$ Hence, $k_{0}=\underline{k}=0.5$ and $k_{N_{k}}=%
\overline{k}=1.5$.

By solving the Lippmann-Schwinger equation (\ref{2.192}) for every point $%
k_{m}$ and taking the values of $u(0,k_{m})$ we obtain the noiseless
boundary function $g_{0}(k_{m})$ in (\ref{2.8}). Next, we add the random
noise in this function 
\begin{equation}
g_{0,\delta }(k_{m})=g_{0}(k_{m})(1.0+\delta \sigma (k_{m})),\,\sigma
=\sigma _{r}(k_{m})+i\sigma _{r}(k_{m}),\,m=0,\dots N_{k},  \label{3}
\end{equation}%
where $\delta $ is the noise level, $\sigma _{r}$ and $\sigma _{i}$ are
random floating point numbers, uniformly distributed between $-1.0$ and $1.0$%
. In our computations we use $\delta =0.05$, i.e. our data have $5\%$ of
noise. To reduce this noise the function $g_{0,\delta }(k_{n})$ is smoothed
out by using the standard averaging procedure. Next, we calculate the
boundary function $q_{0,\delta }$ and its derivative $q_{1,\delta }$ (\ref%
{eq:v_bcs}).

%
%
%

\subsection{Location estimation and data propagation}

\label{sec:6.2}

We have discovered in our computations that prior the minimization of the
functional (\ref{3.14}), it is important to estimate the location of the
target first. The same observation was made in \cite{KlibKol2} for the
3-D case: see Figure 2 in \cite{KlibKol2}.
Procedures of \cite{KlibKol2} for these estimates are different from
the one described here, both of them are heuristic ones. The procedure
described below is also a \emph{heuristic} one.

First, we recall one of steps of \cite{KlibKol1}. Let $v\left( x,k\right) $
be the function defined in (\ref{eq:v}). Rather than considering the
truncated Fourier series (\ref{3.7}), one can denote in (\ref{3.6}) $s\left(
x,k\right) =v_{k}\left( x,k\right) .$ Hence, 
\begin{equation}
v\left( x,k\right) =-\displaystyle\int\limits_{k}^{\overline{k}}s\left(
x,\tau \right) d\tau +v(x,\overline{k}).  \label{100}
\end{equation}%
Then substitution of (\ref{100}) in equation (\ref{3.6}) results in a
nonlinear integro differential equation with respect to the function $%
s\left( x,k\right) .$ The function $v(x,\overline{k})$ is also unknown in
this equation. In \cite{KlibKol1} this function is called the
\textquotedblleft tail function". It follows from (\ref{2.22})-(\ref{2.24})
that 
\begin{equation}
v(x,k)=\frac{r\left( x\right) }{k}+O\left( \frac{\ln \left( 2k\right) }{k^{2}%
}\right) ,k\rightarrow \infty .  \label{101}
\end{equation}%
Assuming that the number $\overline{k}>>1$ is sufficiently large, dropping
the second term in the right hand side of (\ref{101}), substituting $v(x,%
\overline{k})=r\left( x\right) /\overline{k}$ in that integro differential
equation and setting in it $k:=$ $\overline{k},$ one obtains that $r^{\prime
\prime }\left( x\right) =0.$ Solution of the latter equation for $x\in
\left( 0,1\right) $ with certain boundary conditions at $x=0,1$ provides an
approximation for the tail function.

\begin{figure}[tbp]
\begin{center}
\includegraphics[width=0.33\textwidth]{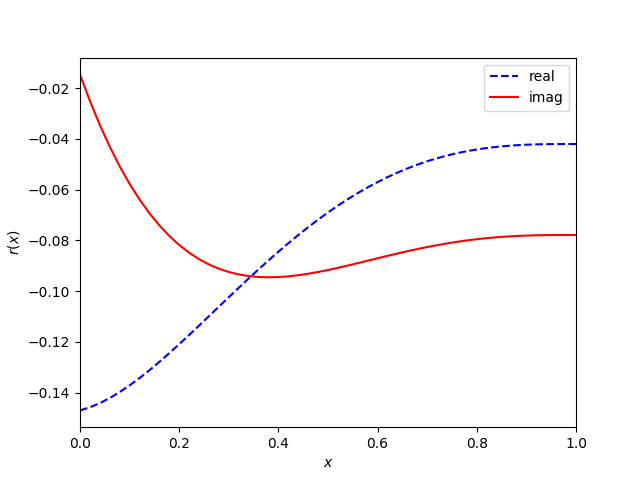}
\end{center}
\caption{The real (dashed line) and imaginary (solid line) parts of the function 
$r(x)$ for the target with $\protect\widehat{c}_{true}=5.0$ and $x_{loc}=0.4$%
}
\label{fig:tail}
\end{figure}

Unlike \cite{KlibKol1}, to estimate the location of the target from our
computationally simulated data, we do not assume that $\overline{k}$ is
sufficiently large. Nevertheless, we solve the same equation $r^{\prime
\prime }\left( x\right) =0$. More precisely, we solve the following boundary
value problem: 
\begin{equation}
r^{\prime \prime }=0,\quad x\in (0,1),  \label{4}
\end{equation}%
\begin{equation}
r(0)=q_{0,\delta }(\overline{k}),\quad r^{\prime }(0)=q_{1,\delta }(%
\overline{k}),\quad r^{\prime }(1)=0.  \label{5}
\end{equation}%
Since the problem (\ref{4}), (\ref{5}) is overdetermined, we solve it by the
quasi-reversibility method (QRM). We refer to \cite{KlibLoc} for details
about this specific version of QRM. More precisely, we minimize the
following functional $I_{\alpha }(r)$: 
\begin{equation}
I_{\gamma }(r)=\frac{1}{2}\left( \Vert r^{\prime \prime }\Vert
_{L^{2}(0,1)}+\gamma \Vert r\Vert _{L^{2}(0,1)}\right) ,  \label{eq:loc_est}
\end{equation}%
subject to boundary conditions (\ref{5}). Here $\gamma >0$ is the
regularization parameter. Figure \ref{fig:tail} depicts the so computed
function $r(x)$ for the target with $\widehat{c}_{true}=5.0$ and $%
x_{loc}=0.4 $, see (\ref{eq:c_true}). Note that the minimal value of $\func{%
Im}r\left( x\right) $ is near the point $x=0.4,$ which is the center $%
x_{loc} $ of this simulated target. Therefore, we have estimated the
location of the center of this target $x_{est}$. We have seen in all 16 our
computational experiments that this approach always works: it gives us the
value of $x_{est}$ with an error of about $5\%$.

Next, in cases when the estimated location $x_{est}>0.1$, we propagate the
data function $g_{0,\delta }\left( k\right) $ from the point $x=0.0,$ where
this function is given, to the point $x_{tar}=x_{est}-0.1.$ By (\ref%
{eq:c_true}) this means that the maximal \textquotedblleft allowable" width
of the target in our case is 0.2.

We point out that a 3-D version of the data propagation procedure has been
widely used in the previous works of Klibanov with coauthors for both
computationally simulated and experimental data, see, e.g. \cite%
{KlibKol2,Kfreq} and references cited therein. Since $c(x)=1.0$ for $%
x\in (0,x_{tar})$ and $x_{0}<0,$ then by (\ref{2.4}) $u^{\prime \prime
}+k^{2}u=0$ for $x\in (0,x_{tar}).$ Hence, 
\[
u\left( x,k\right) =D_{1}\left( k\right) e^{ikx}+D_{2}\left( k\right)
e^{-ikx},\, x\in \left( 0,x_{tar}\right) , 
\]%
where complex numbers $D_{1}\left( k\right) ,D_{2}\left( k\right) $ depend
only on $k$. To find these numbers, we need to know $u\left( 0,k\right) $
and $u^{\prime }\left( 0,k\right) .$ The latter numbers, in turn can be
easily found from formulae (\ref{2.60})-(\ref{2.160}). Hence, we perform the
data propagation via the following formulae:

\[
g_{prop}(x_{tar},k)=\frac{u_{p}(x_{tar},k)}{u_{0}(x_{tar},k)},\quad
u_{prop}(x_{tar},k)=D_{1}e^{ikx_{tar}}+D_{2}e^{-ikx_{tar}}, 
\]%
\[
D_{1}=u_{0}(0)\left( g_{0}-1.0\right) ,\quad D_{2}=u_{0}(0), 
\]%
where $k-$dependent functions $g_{prop}(x_{tar},k)$ and $u_{prop}(x_{tar},k)$
play now the role of functions $g_{0}\left( k\right) $ and $u\left(
0,k\right) $ in which $\left\{ x=0\right\} $ is replaced with $\left\{
x=x_{tar}\right\} .$

\subsection{The optimal number $N$ terms in the expansion (\protect\ref{3.7})%
}

\label{sec:6.3}

\begin{figure}[tbp]
\begin{center}
\subfloat[\label{fig:c_basis}]{\includegraphics[width=0.33%
\textwidth]{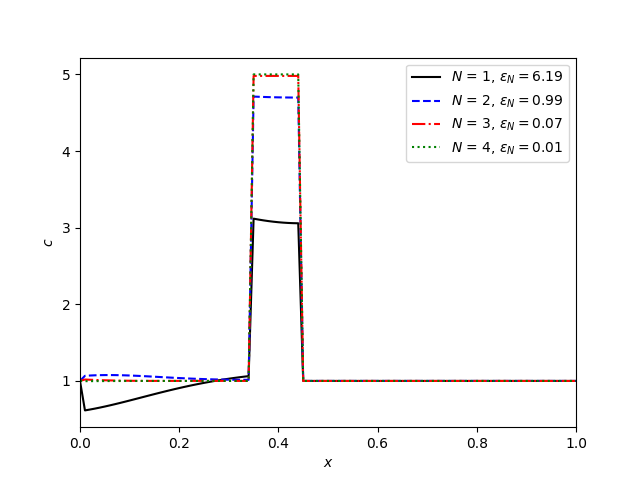}} \subfloat[\label{fig:basis}]{%
\includegraphics[width=0.33\textwidth]{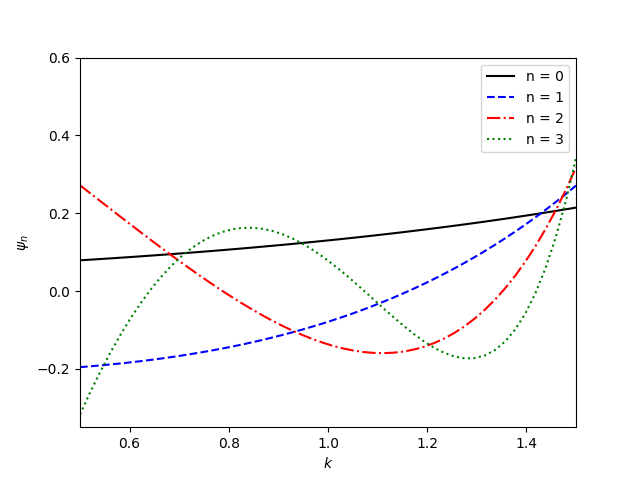}}
\end{center}
\caption{a) The approximate functions $c_{appr,N}(x)$, b) basis functions $%
\protect\phi_n(k)$}
\label{fig:2}
\end{figure}

We need to determine the optimal number $N$ of terms in the truncated
Fourier series (\ref{3.7}). To do this, we first solve the
Lippmann-Schwinger equation (\ref{2.192}) for a reference target with $%
c\left( x\right) :=\widehat{c}_{true}\left( x\right) ,$ see (\ref{eq:c_true}%
) for $\widehat{c}_{true}.$ This way we obtain the functions $w_{true}$ in (%
\ref{2.100}) and $v_{true}(x,k)$ in (\ref{eq:v}). Next, we compute vector
functions $y_{true,N}(x)$ in (\ref{3.7}), (\ref{3.90}) for different values
of $N$ and reconstruct approximate functions $c_{appr,N}(x)$ via (\ref{3.8}%
). Figure \ref{fig:c_basis} shows the functions $c_{appr,N}(x)$ for $N=$ $%
1,2,3,4$, where $\widehat{c}_{true}=5.0$ and $x_{loc}=0.4$. The
corresponding basis functions $\psi _{n}(k)$ are shown on Figure \ref%
{fig:basis}. One can see that the functions $c_{true,N}(x)$ are accurately
approximated for both $N=3$ and $N=4$, and their approximation errors $%
\varepsilon _{N}=\Vert c_{appr,N}-c_{true}\Vert _{L^{2}(0,1)}$ are
sufficiently small: $\varepsilon _{N}=0.07$ and $0.01$, respectively.
Therefore, we choose the optimal number of functions in our basis $N=3$.

\subsection{Numerical implementation}

\label{sec:6.4}

For the numerical solution we use the finite difference discretization
method. So we divide the interval $x\in \lbrack 0,1]$ into $N_{x}$ equal
subintervals and obtain the mesh $x_{j}=jh_{x}$, $j=0,\dots ,N_{x}$, $%
h_{x}=1.0/N_{x}$ . Combining this mesh with the mesh for wave numbers $k_{m}$
defined in section 6.1, we obtain the $N_{k}\times N_{x}$ two dimensional
mesh with the grid points $(k_{m},x_{j})$. We need to find the discrete
function $V=\{v_{m,j}\}$, where $v_{m,j}=v(k_{m},x_{j})$ are its values at
those grid points. The discrete version of equation (\ref{3.7}) can be
written as 
\begin{equation}
V=\Psi Y,  \label{eq:V}
\end{equation}%
where $\Psi =\{\psi _{m,n}\}$ is the $N_{k}\times N$ matrix with $\psi
_{m,n}=$ $\psi _{n}(k_{m})$ and $Y=\{y_{n,j}\}$ is the $N\times N_{x}$ two
dimensional discrete vector function $y$ (\ref{3.90}). 

The main objective of the proposed method is to find the minimizer of
functional $J_{\lambda ,\alpha }(y)$ (\ref{3.14}). Naturally, \ we minimize this functional
numerically for the case of the discrete function $Y$. The details of finite
difference discretization of the functional and its gradient are not
described here for brevity. We refer to \cite{Kuzh}, where discretization of
similar functional is performed analytically using the Kronecker delta
function. Although, our theory suggest the gradient projection method for
the minimization, we have noticed that the conjugate gradient method (CG)
works well for our problem. This method is easier to implement numerically
and it gives practically the same results as the gradient projection method.
The latter has been consistently observed in the previous works on the numerical issues of the convexification method \cite{KlibKol1,  KlibThanh, KlibKol2, KlibYag, KEIT, Bak}. Also, to decrease the
necessity of calculating the functional and its gradient on each iterative
step, we have decided not to use the standard line search algorithm to seek
the step size of the minimization process. Instead, we begin with the
initial step size $10^{-7}$. On each iteration, the step size is reduced $10$
times for the next iteration if the value of the functional on the current
iteration is greater than its value on the previous iteration. Otherwise,
the step size remains the same for the next iteration. Also, after every
1000 iterations the step size is increased by 10 times. The minimization
algorithm is stopped either after 15000 iterations or when the step size
becomes less then $10^{-14}$. The latter means that the functional can no
longer decrease and its minimizer is reached.

\subsection{Algorithm}

\label{sec:6.5}

Here, we summarize our algorithm for reconstructing the unknown function $%
c(x)$ from noisy computationally simulated data $g_{0,\delta }\left(
k\right) $. Below in this section, we assume that our data have noise and
omit the subscript $\delta $ for brevity. Also, all functions in this
section are discrete, unless otherwise specified.

\begin{enumerate}
\item Estimate the location of target and propagate data if necessary, i.e.
if $x_{est}>0.1$.

\item Calculate the boundary conditions $q_{0}$, $q_{1}$ in (\ref{eq:v_bcs})
and then $f_{0}$, $f_{1}$ in (\ref{3.12}) as 
\[
{f_{0,n}}=\sum_{m=0}^{N_{k}}\psi _{m,n}{q_{0,m}},\quad {f_{1,n}}%
=\sum_{m=0}^{N_{k}}\psi _{m,n}q_{1,m},n=0,...,N-1, 
\]%
\[
f_{0}=\left( f_{0,0},...,f_{0,N-1}\right) ^{T},\, f_{1}=\left(
f_{1,0},...,f_{1,N-1}\right) ^{T}. 
\]

\item Define the initial guess $Y_{0}$ for the vector function $Y$ as $%
Y_{0}=f$ in (\ref{2}).

\item Minimize the functional $J_{\lambda ,\alpha }(Y)$ in (\ref{3.14}).
Then transform the found vector function $Y$ in the vector function $V$ in (%
\ref{eq:V}).

\item Using (\ref{3.3}) with $k=\underline{k},$ compute the approximation $%
\beta _{comp} $.

\item After averaging $\beta _{comp}\left( x\right) ,$ determine the
coefficient $c_{comp}$ as follows: 
\begin{equation}
c_{comp}=\left\{ 
\begin{array}{cc}
\mbox{Re}(\beta _{comp})+1.0, & \mbox{if }\mbox{Re}(\beta _{comp})\geq \rho
\max (\mbox{Re}(\beta _{comp})), \\ 
1.0, & \mbox{otherwise}%
\end{array}%
\right.  \label{eq:c_comp}
\end{equation}%
where $\rho \in (0,1)$ is the truncated factor.
\end{enumerate}

Step 6 of the algorithm is a simple 
post-processing to reconstruct our simulated targets.

\subsection{Reconstruction results for computationally simulated data}

\label{sec:6.6}

\begin{figure}[tbp]
\begin{center}
\subfloat[\label{fig:3_1}]{\includegraphics[width=0.33\textwidth]{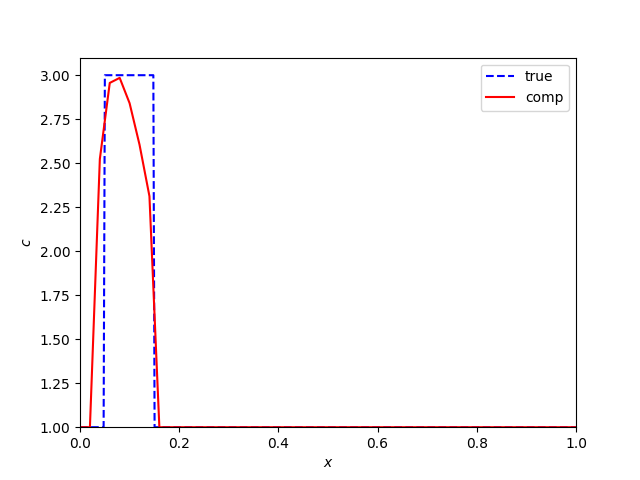}} %
\subfloat[\label{fig:4_2}]{\includegraphics[width=0.33\textwidth]{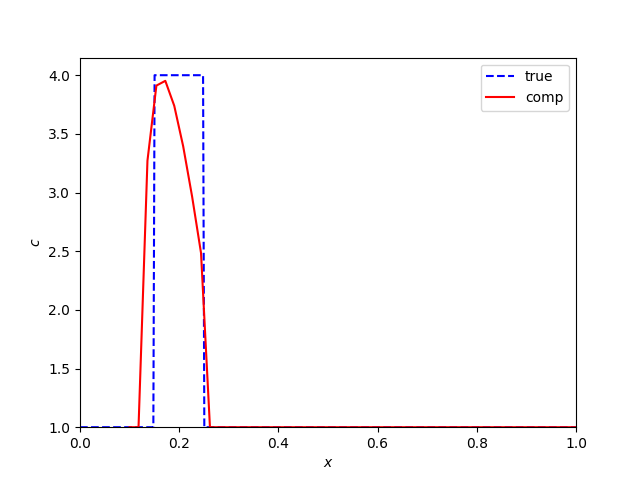}} \\[%
0pt]
\subfloat[\label{fig:5_3}]{\includegraphics[width=0.33\textwidth]{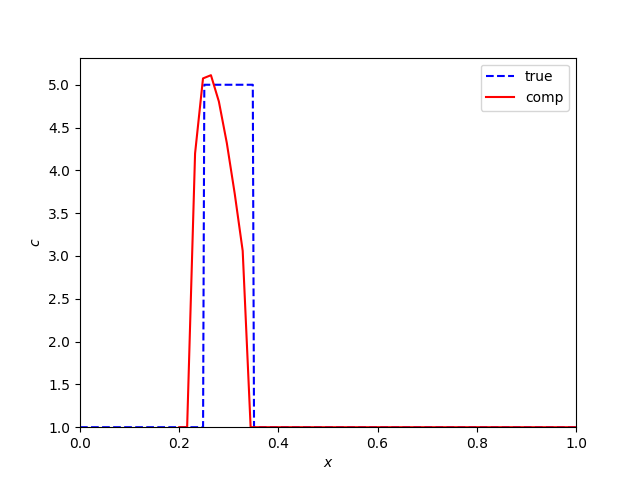}} %
\subfloat[\label{fig:6_4}]{\includegraphics[width=0.33\textwidth]{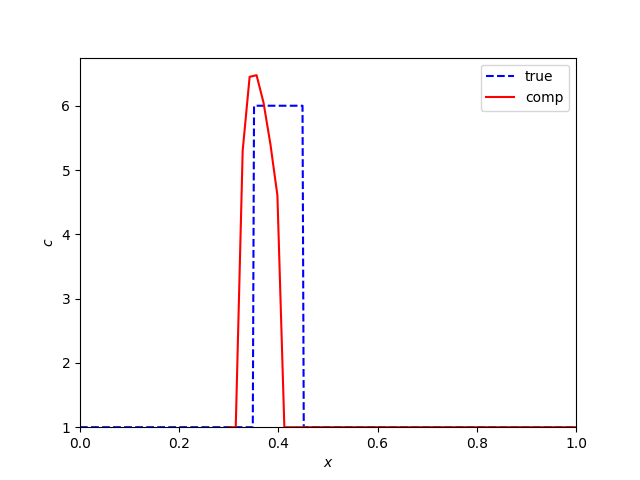}}
\end{center}
\caption{Reconstruction results for targets with: a) $\protect\widehat{c}%
_{true} = 3.0$, $x_{loc}=0.1$, b) $\protect\widehat{c}_{true} = 4.0$, $%
x_{loc}=0.2$, c) $\protect\widehat{c}_{true} = 5.0$, $x_{loc}=0.3$, d) $%
\protect\widehat{c}_{true} = 6.0$, $x_{loc}=0.4$}
\label{fig:res_sim}
\end{figure}

In this section we present the results of reconstructions via the proposed
algorithm for the computationally simulated targets with $\widehat{c}_{true}$
and $x_{loc}$ set as in (\ref{eq:targers}). These results are obtained using
the Carleman weight function $e^{-2\lambda x}$ with the parameter $\lambda
=3.0$ in (\ref{3.14}), the regularization parameter $\alpha =0.05$, the
discretization parameters $N_{x}=50$, $N_{k}=3$, $N=3$, and the truncation
factor $\rho =0.5$ in (\ref{eq:c_comp}). Table \ref{tab:res_sim} lists all
reconstructed targets with the maximum value of the computed function, 
\begin{equation}
\widehat{c}_{comp}=\max (c_{comp}).  \label{6.1}
\end{equation}
We define the relative computational error: 
\[
\varepsilon _{comp}=\frac{|\widehat{c}_{comp}-\widehat{c}_{true}|}{\widehat{c%
}_{true}}\cdot 100\%. 
\]%
The reconstructed functions $c_{comp}$ in (\ref{eq:c_comp}) for some
selected targets are shown in Figure \ref{fig:res_sim}. We can see that
these results are quite accurate ones, given that with the noise level in
the data is 5\%.


\begin{table}[!h]
\caption{Reconstruction results for simulated targets.}
\label{tab:res_sim}
\begin{center}
\begin{tabular}{|c|c|c|c||c|c|c|c|}
\hline
$\widehat{c}_{true} $ & $x_{loc}$ & $\widehat{c}_{comp}$ & $%
\varepsilon_{comp}, \%$ & $\widehat{c}_{true} $ & $x_{loc}$ & $\widehat{c}%
_{comp}$ & $\varepsilon_{comp}, \%$ \\ \hline
3.0 & 0.1 & 2.98 & 0.67 & 5.0 & 0.1 & 5.32 & 6.40 \\ 
& 0.2 & 3.13 & 4.33 &  & 0.2 & 5.14 & 2.80 \\ 
& 0.3 & 2.80 & 6.67 &  & 0.3 & 5.11 & 2.20 \\ 
& 0.4 & 3.17 & 5.67 &  & 0.4 & 5.19 & 3.80 \\ \hline
4.0 & 0.1 & 4.28 & 7.00 & 6.0 & 0.1 & 6.19 & 3.17 \\ 
& 0.2 & 3.95 & 1.25 &  & 0.2 & 6.25 & 4.17 \\ 
& 0.3 & 4.03 & 0.75 &  & 0.3 & 6.39 & 6.50 \\ 
& 0.4 & 4.12 & 3.00 &  & 0.4 & 6.47 & 7.83 \\ \hline
\end{tabular}%
\end{center}
\end{table}

\subsection{Reconstruction results for experimental data}

\label{sec:6.7}

Now we demonstrate the reconstruction results for the experimental data.
Recall that these data were collected in the field (as opposed to a
laboratory) by the Forward Looking Radar of the US Army Research Laboratory 
\cite{Radar}. The scheme of data collection is presented on Figure \ref%
{fig:setup}. Originally the backscattering time dependent data are measured,
one time resolved curve for each target. To obtain the data in the frequency
domain, we apply Fourier transform. Since samples of shapes of both time
dependent experimental curves and their Fourier transforms can be found in 
\cite{KlibLoc}, we do not display them here.

The same experimental data were used in our previous works \cite{KlibKol1, KlibLoc, Karch, Kuzh2, IEEE}. Hence, we do not describe here details
of the data collection method as well as the data pre-processing procedure.
Instead, we refer to the above mentioned previous works. Especially detailed
description can be found in \cite{Karch,Kuzh2}. In references \cite{KlibLoc, Karch, Kuzh2, IEEE} these data were treated by the tail
functions globally convergent method and in \cite{KlibKol1} they were
treated by the previous version of the convexification method.

Targets of our interest were surrounded by cluttered environment, which is
of course a complicating factor for their imaging. Horizontal coordinates of
targets were provided by Ground Positioning System (GPS) with a good
accuracy. As to the burial depths of targets, they are not of an interest in
this specific application. In addition, we had the data for two targets
located in air. Hence, the burial depth for these two makes no sense.
Furthermore, it is clear from the descriptions of \cite{Karch,Kuzh2} of the
data collection process that it is unlikely that an information about burial
depths of buried targets can ever be extracted from these data. All what was
known to us was that burial depths of targets buried in the ground (dry
sand) was a few centimeters. With respect to our algorithm this means that
we are not interested in the location of the target. Rather, our interest is
in the target/background contrast in the dielectric constant, see (\ref{6.2}%
). Hence, unlike the computationally simulated data, we do not apply here
the procedure of section 6.2.

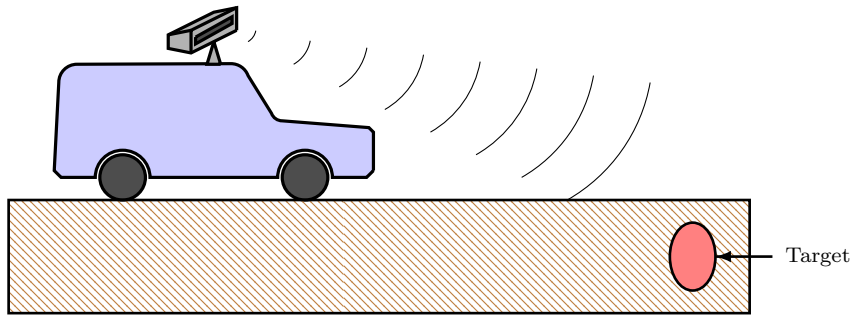
\begin{figure}[tbp]
\begin{center}
\begin{tikzpicture}[font=\footnotesize,scale=0.3]  
\draw [fill=blue!20!white,line width=1.2]  (11,0.3) -- (11,2) -- (10.8,2.2) -- (7,2.5) arc (270:200:0.5) -- (5.5,4.5) arc (20:90:0.8)  -- (-2,5) arc (90:180:0.8) -- (-3,0.3) --  (-2.7,0) -- (-1.2,0) arc (180:0:1.2) -- (6.8,0) arc (180:0:1.2) -- (10.7,0) -- (11,0.3) ;
\draw [fill=black!70!white,line width=1.2] (0,0) circle (1);
\draw [fill=black!70!white,line width=1.2] (8,0) circle (1);

\draw [fill=black!30!white,line width=1.2] (4.3,5) -- (4,6) -- (3.7,5) -- (4.3,5);

\draw [fill=black!30!white,line width=1.2] (5,6.5) -- (5,7.5) -- (4,7.3) -- (2, 6.3) -- (2, 5.7) --   (3,5.5) -- (5,6.5);

\draw [line width=1.2] (3,5.5) -- (3,6.5) -- (5,7.5);
\draw [line width=1.2] (3,6.5) -- (2,6.3);

\draw [fill=black!70!white,line width=1.2] (3.2,6.2) -- (4.8,7.1) -- (4.8,6.8) -- (3.2,5.9) -- (3.2,6.2);

\draw [pattern=north west lines, pattern color=brown,line width=1.0] (-5,-1) rectangle (27.5,-6);

\draw[fill=red!50!white,line width=1.0]   (25,-3.5) ellipse (1 and 1.5);


\node  at (30.5,-3.5) {Target}; 

\draw[onearrow={latex},line width=1.0] (28.5,-3.5) -- (26,-3.5); 

\draw (5.5,6) arc (300:350:0.7);
\draw (7.5,5.0) arc (300:350:1.5);
\draw (9.5,4) arc (300:350:2.5);
\draw (11.5,3) arc (300:350:3.5);
\draw (13.5,2) arc (300:350:4.5);
\draw (15.5,1) arc (300:350:5.5);
\draw (17.5,0) arc (300:350:6.5);
\draw (19.5,-1) arc (300:350:7.5);
\end{tikzpicture}
\end{center}
\caption{Data collection scheme}
\label{fig:setup}
\end{figure}

We posses experimental data for five (5) targets. Our \emph{a priori}
information was that two (2) targets, bush and a wood stake, were located in
air, and three (3) targets, metallic box, metallic cylinder and plastic
cylinder, were buried in a dry sand. We also knew that each experimental
data set was collected for a single target only. We introduce the number $%
c_{bg}$ for the dielectric constant of the background, which is $c_{bg}=1.0$
if the target is in air, and $c_{bg}\in \left[ 3.0,5.0\right] $ if the
target is buried. Here we use values of the dielectric constant of dry sand, see tables \cite{Tables} of dielectric constants.
Consequently, the computed number $\widehat{c}_{comp}=\max c_{comp}(x)$ in (%
\ref{6.1}) is an estimation of the contrast between the target and the
background $c_{contrast}$, 
\begin{equation}
c_{contrast}=\frac{c_{target}}{c_{bg}}\approx \widehat{c}_{comp}.
\label{6.2}
\end{equation}%
where $c_{target}$ is the true dielectric constant of the target and the
function $c_{comp}(x)$ is as in (\ref{eq:c_comp}).

Also, we know that the dielectric constant of plastic cylinder is less then
the one of sand with $0<c_{contrast}<1.0$, see \textquotedblleft plastic
pellets" in \cite{Tables}. Therefore, we modify (\ref{eq:c_comp}) for this
case as: 
\[
c_{contrast}=\left\{ 
\begin{array}{cc}
\mbox{Re}(\beta _{comp})+1.0, & \mbox{if }\mbox{Re}(\beta _{comp})\leq \rho
\min (\mbox{Re}(\beta _{comp})), \\ 
1.0, & \mbox{otherwise and also if }\func{Re}\left( \beta _{comp}\right)
<-1.0,%
\end{array}%
\right. 
\]%
where $\rho = 0.5$.  Here, the condition $\func{Re}(\beta _{comp})<-1.0$
must be applied before the truncation. Next, we introduce the number $%
\widehat{c}_{est}$, which is our estimation of the dielectric constant of
the target, 
\begin{equation}
\widehat{c}_{est}=c_{bg}\widehat{c}_{comp},  \label{6.3}
\end{equation}%
where $\widehat{c}_{comp}=\min (c_{comp})$ for the plastic cylinder, and $%
\widehat{c}_{comp}=\max (c_{comp})$ for other targets. \emph{A priori}
differentiation between the plastic cylinder and other targets is performed
prior our computations via an analysis of the data, see \cite{Kuzh2}.

The numbers function $\widehat{c}_{comp}$ as well as estimates of dielectric
constant $\widehat{c}_{est}$ via (\ref{6.3}) for our targets are listed in
Table \ref{tab:res_exp}. Since we did not measure dielectric constants of
targets, then the values of $c_{target}$ were obtained from tables \cite{Tables}. Also, it was shown in \cite{Kuzh2} that one can assign large
values of dielectric constants to metallic targets as $c_{met}:=c_{target}%
\in \left[ 10,30\right] $. As to the bush, this target is considered as the
hardest one to image since it is obviously a very heterogeneous one due to
the presence of many leaves. We took for this case $c_{target}$ from \cite%
{Chuah}. Our results are summarized in Table \ref{tab:res_exp}.

\begin{figure}[tbp]
\begin{center}
\subfloat[\label{fig:bush}]{\includegraphics[width=0.33\textwidth]{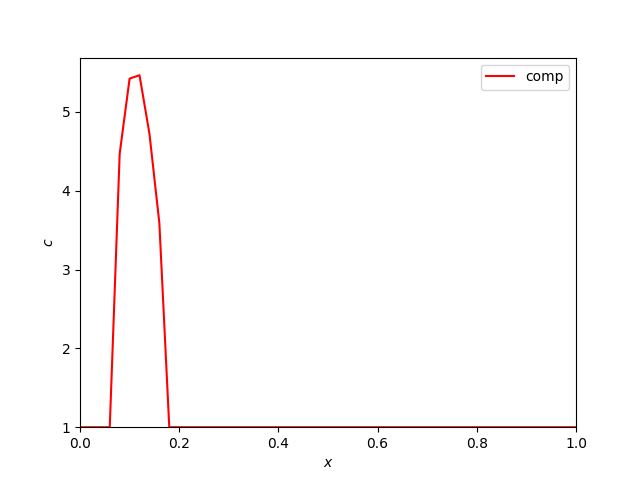}} %
\subfloat[\label{fig:wood}]{\includegraphics[width=0.33\textwidth]{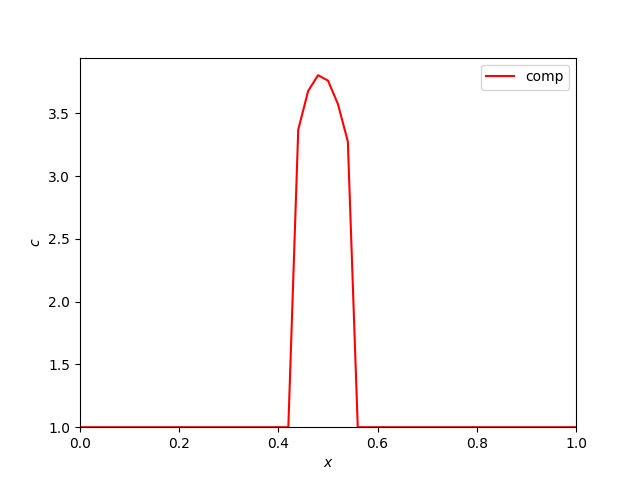}} 
\\[0pt]
\subfloat[\label{fig:mbox}]{\includegraphics[width=0.33\textwidth]{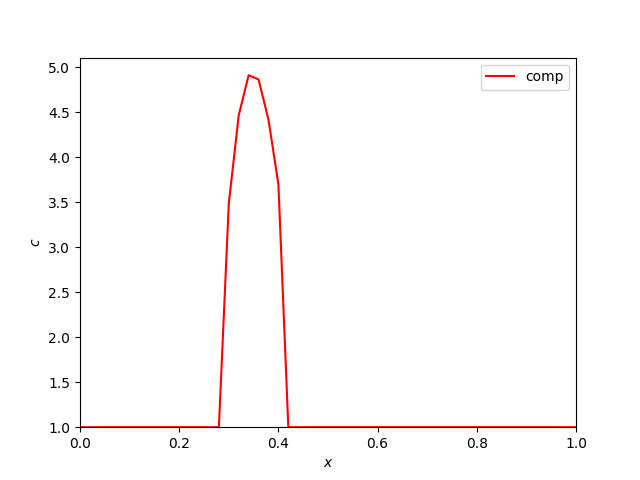}} %
\subfloat[\label{fig:plastic}]{\includegraphics[width=0.33%
\textwidth]{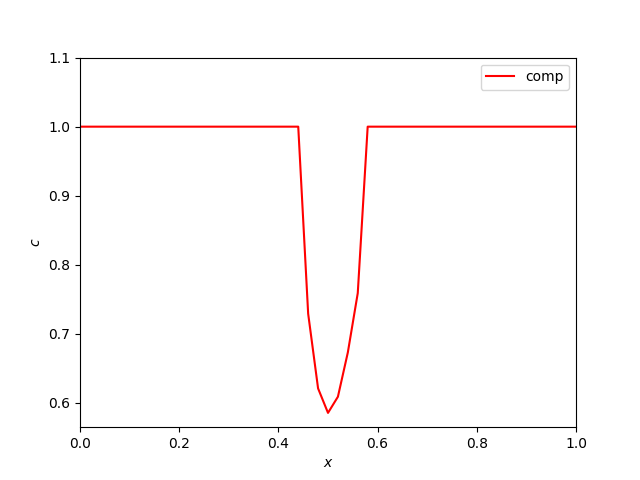}}
\end{center}
\caption{Reconstruction results for: a) bush, b) wood stake, c) metallic box, d)
plastic cylinder}
\label{fig:res_exp}
\end{figure}

\begin{table}[h]
\caption{Reconstruction results for experimental targets.}
\label{tab:res_exp}
\begin{center}
\begin{tabular}{|l|c|c|c|c|c|}
\hline
Target & Air/Sand & $\widehat{c}_{contrast}$ & $c_{bg}$ & $\widehat{c}_{est}$
& $c_{target}$ \\ \hline
Bush & Air & 5.47 & 1.0 & 5.47 & $[3.0,20.0]$ \\ 
Wood stake & Air & 3.80 & 1.0 & 3.80 & $[2.0,6.0]$ \\ 
Metal box & Sand & 4.91 & $\left[ 3.0,5.0\right] $ & $\left[ 14.73, 24.55%
\right]$ & $[10.0,30.0]$ \\ 
Metal cylinder & Sand & 4.84 & $\left[ 3.0,5.0\right] $ & $\left[ 14.52,
24.20\right] $ & $[10.0,30.0]$ \\ 
Plastic cylinder & Sand & 0.59 & $\left[ 3.0,5.0\right] $ & $\left[ 1.77,
2.95\right] $ & $[1.1,3.2]$ \\ \hline
\end{tabular}%
\end{center}
\end{table}

One can see that all values of estimated dielectric constants $\widehat{c}%
_{est}$ in this table are within tabulated limits. Given that dielectric
constants of targets were not measured in experiments, Table \ref%
{tab:res_exp} is a quite encouraging one for the engineering part of this
research group (AS and LN). Indeed, results presented in this table indicate
that a future software based on the above algorithm might indeed provide
rather accurate estimates of dielectric constants of targets of interest. An
important additional point of the encouragement of engineers is that these
results are obtained for targets surrounded by a realistic cluttered
environment. Hence, engineers conjecture that an intriguing opportunity
might occur in the future to decrease the false alarm rate, as mentioned in
section 1. Certainly more comprehensive studies of large collections of
experimental data sets are necessary to verify this conjecture.

\section*{Acknowledgements}

The work of M. V. Klibanov and A.E. Kolesov was supported by the US Army
Research Laboratory and US Army Research Office grant W911NF-15-1-0233 as
well as by the Office of Naval Research grant N00014-15-1-2330.

%
%

\end{document}